\documentclass[12pt]{article}

\usepackage{amssymb,amsmath,amsthm,amsfonts}
\usepackage{graphicx,graphics}
\usepackage{epstopdf}
\usepackage{algorithm,algpseudocode,amsopn}
\graphicspath{{./pics/}}
\usepackage{verbatim}
\usepackage{mathtools}
\usepackage{threeparttable}
\usepackage{booktabs}
\usepackage{url}
\usepackage[nohead,margin=1.1in]{geometry}
\usepackage{color}
\usepackage{tikz}
\usetikzlibrary{calc}
\numberwithin{equation}{section}
\numberwithin{algorithm}{section}

\newtheorem{theorem}{Theorem}[section]
\newtheorem{lemma}{Lemma}[section]
\newtheorem{proposition}{Proposition}[section]
\newtheorem{corollary}{Corollary}[section]
\newtheorem{definition}{Definition}[section]

\allowdisplaybreaks

 \newcommand{\ind}{\,\mbox{d}}

\title{On the recovery of two function-valued coefficients in the Helmholtz equation for inverse scattering problems via inverse Born series\thanks{The research of F.C. is partially supported by NSF Grant DMS 2406313. }}

\author{Fioralba Cakoni\footnotemark[2]\thanks{Department of Mathematics,
Rutgers University, 110 Frelinghuysen Road, Piscataway, NJ 08854-8019,  USA. (\texttt{fc292@math.rutgers.edu}, \texttt{zz569@math.rutgers.edu})} \and Shixu Meng\footnotemark[3]\thanks{Department of Mathematics, Virginia Tech, 24061 Blacksburg,  USA.  (\texttt{sgl22@vt.edu})} \and Zehui Zhou\footnotemark[2]}

\date{}

\begin{document}

\maketitle
\begin{abstract}
In this work, we construct the Born and inverse Born approximation and series to recover two function-valued coefficients in the Helmholtz equation for inverse scattering problems from the scattering data at two different frequencies. 
An analysis of the convergence and approximation error of the proposed regularized inverse Born series is provided. 
The results show that the proposed series converges when the inverse Born approximations of the perturbations are sufficiently small. 
The preliminary numerical results show the capability of the proposed regularized inverse Born approximation and series for recovering the isotropic inhomogeneous media.

\end{abstract}

\section{Introduction}\label{sec:intro}

The scattering of the time-harmonic incident plane wave $u^i(x;\hat{\theta};k):=e^{ik x\cdot \hat{\theta}}$
with the probing wave frequency $k$ and incident direction $\hat{\theta}\in\mathbb{S}$ (with $\mathbb{S}$ denoting the unit sphere) by an inhomogeneity occupying a bounded Lipschitz domain $\Omega \subset \mathbb{R}^2$ is mathematically formulated {as follows}: find the total wave field $u\in H_{loc}^1(\mathbb{R}^2)$ with $u=u^s+u^i$ such that
\begin{align}
&\nabla\cdot a(x) \nabla u+k^2 n(x) u=0  \;\;\mbox{in}\;\; \mathbb{R}^2,\label{eqn:prob}\\
&\lim_{r\to \infty} r^\frac12 \left(\dfrac{\partial u^s}{\partial r}-ik u^s\right)=0 \;\;\mbox{with}\;\;r=|x|  \label{eqn:src}
\end{align}
where 
\begin{align*}
a(x)=&\left\{
\begin{array}{cl}
    1 & \mbox{in }\; \mathbb{R}^2\setminus \overline{\Omega}\\
    1+\gamma(x) & \mbox{in }\; \overline{\Omega}
\end{array}
\right. \quad \mbox{ and }\quad 
n(x)=\left\{
\begin{array}{cl}
    1 & \mbox{in }\; \mathbb{R}^2\setminus \overline{\Omega}\\
    1+\eta(x) & \mbox{in }\; \overline{\Omega}
\end{array}
\right.,
\end{align*}
with the contrasts $\gamma$ being piecewise continuously differentiable and $\eta$ being bounded in $\overline\Omega$, 
and the Sommerfeld radiation condition \eqref{eqn:src} on the scattered wave field $u^s$ is satisfied uniformly with respect to $\hat{x}=x/|x|$. 
{Problem \eqref{eqn:prob}--\eqref{eqn:src} models} the scattering of acoustic waves by an isotropic inhomogeneous medium with contrasts in both the sound speed and density \cite{colton-kress} in two dimensions. 
Thanks to Sommerfeld radiation condition \eqref{eqn:src}, the scattered filed $u^s$ assumes the following asymptotic behavior
$$u^s(x;\hat{\theta};k)=\frac{e^{ik |x|}}{\sqrt{|x|}}u^\infty(\hat{x},\hat{\theta};k)+\mathcal{O}\left(\frac{1}{|x|^{3/2}}\right),$$
where $u^\infty(\hat{x},\hat{\theta};k)$, known as the far-field pattern of the scattered field, is a function of $\hat x \in {\mathbb S}$ for fixed $k$ and $\hat{\theta}$.
\begin{definition}[Scattering Data and the Inverse Problem] {\em The set of measured far-field patterns
$$\left \{u^\infty(\hat{x},\hat{\theta};k): \, \hat x\in {\mathbb S}, \hat{\theta}\in{\mathbb S}\right\}$$
is called the {\it scattering data} at frequency $k$. The {\it inverse scattering problem} we consider here is to determine the contrasts $\gamma$ and $\eta$ from the scattering data at two different frequencies $k_1\neq k_2$.}
\end{definition}

It is well-known that the uniqueness of the inverse scattering problem derived from \eqref{eqn:prob}--\eqref{eqn:src} holds \cite{uniq2}. We also refer the reader to \cite{uni2, uni1} for coefficients less regular than assumed here. 
By noting that the incident wave field  $u^i(x;\hat{\theta};k)=e^{ik x\cdot \hat{\theta}}$ solves the Helmholtz equation with the homogeneous media (i.e. $\gamma\equiv0$ and $\eta\equiv0$) in ${\mathbb R}^2$,  the perturbed Helmholtz equation \eqref{eqn:prob} for the scattered field takes the form 
\begin{align}
\nabla\cdot a(x) \nabla u^s+k^2 n(x)u^s=- \nabla\cdot \gamma(x)\nabla u^i-k^2 \eta(x) u^i \;\;\;\mbox{for}\;\; \; x\in \mathbb{R}^2 \label{eqn:prob1}
\end{align}
or equivalently
\begin{equation}
\Delta u^s+k^2 u^s=- \nabla\cdot \gamma(x)\nabla u-k^2 \eta(x) u \;\;\;\mbox{for}\;\;\; x\in  \mathbb{R}^2. \label{eqn:prob2}
\end{equation}
With the Sommerfeld radiation condition \eqref{eqn:src}, it is easy to verify that $u^s$ satisfies the following Lippmann-Schwinger volume integral equation,
\begin{equation}\label{eqn:lsvi}
u^s(x;\hat{\theta};k)= \nabla_x\cdot \int_{\Omega} G^k(x,y)\gamma(y)\nabla_y (u^s+u^i)(y) {\rm d}y+\int_{\Omega} G^k(x,y)k^2 \eta(y) (u^s+u^i)(y){\rm d}y,
\end{equation}
with $G^k(x,y)=i/4 H_0^{(1)}(k |x-y|)$ being the full-space fundamental solution of the Helmholtz equation  in ${\mathbb R}^2$, where $H_0^{(1)}$ is the Hankel function of the first kind of order zero. 
The volume integral equation \eqref{eqn:lsvi} allows $\gamma$ to have jumps in $\overline\Omega$, and we refer the reader to \cite{Costabel:2015} for the sufficient conditions on its solvability. 

Due to the dependence of $u^s$ on $\gamma$ and $\eta$, equation \eqref{eqn:lsvi} is nonlinear with respect to $\gamma$ and $\eta$. 
When the perturbations $\gamma$ and $\eta$ are sufficiently small, the scattered field $u^s$ can be approximated by the linearized version of \eqref{eqn:lsvi}:
\begin{equation}\label{eqn:linear_us}
    u^s(x;\hat{\theta};k)\approx \nabla_x\cdot \int_{\Omega} G^k(x,y)\gamma(y)\nabla_y u^i(y) {\rm d}y+\int_{\Omega} G^k(x,y)k^2 \eta(y)u^i(y){\rm d}y.
\end{equation}
However, when the perturbations grow, the scattered field is no longer well approximated by a simple linearization \eqref{eqn:linear_us}. 
In this case, one may expand the approximation of $u^s$ from \eqref{eqn:linear_us} (e.g., Born approximation) by the Born series with Born approximation being the first-order term. 
A direct inversion of the Born series is typically ill-posed, and the ill-posedness is carried by the inversion of the Born approximation. 
The inverse Born series (IBS) is designed to iteratively reconstruct the potential in a convergent manner and suggests the solution to the inverse problem as an explicitly computable
functional of the scattering data. 
The IBS extends the validity of the inverse Born methods and is widely applied in various inverse scattering problems, including optical tomography \cite{MarkelSullivanSchotland:2003,MoskowSchotland:2008,MoskowSchotland:2009,MachidaSchotland:2015}, electrical impedance tomography \cite{ArridgeMoskowSchotland:2012}, and acoustic and electromagnetic imaging \cite{kilgoreMoskowSchotland：2012,KilgoreMoskowSchotland:2017}.
The main challenge in employing IBS for inverse scattering problems is ensuring the convergence and stability of the series, which has been studied in \cite{MoskowSchotland:2008,ArridgeMoskowSchotland:2012,MachidaSchotland:2015,KilgoreMoskowSchotland:2017,HoskinsSchotland:2022}; see \cite[Chapter 12]{RamlauScherzer:2019} for an overview. 
We also refer the reader to \cite{DesaiLahivaaraMonk:2025} where the applicability of the Born approximation is broadened by two different neural network-based algorithms when the scattering is not
weak.

In this work, we construct the IBS to recover two function-valued coefficients in the Helmholtz equation for inverse scattering problems and provide sufficient conditions on the convergence of the IBS for this problem, when the regularized inverse Born approximation is employed. 
We then discuss the convergence of the IBS with the spectral cutoff regularization based on the disk prolate spheroidal wave functions (PSWFs) as eigenfunctions \cite{meng23data,ZhouAudibertMengZhang:2024}. 
Moreover, when the forward and inverse Born series converge, we derive an upper bound on the approximation error of the IBS.
It is worth noting that the case we consider in this work, involving two unknown function-valued coefficients, is fundamentally distinct from the single-coefficient scenario. 
This distinction introduces challenges for both theoretical and numerical analyses.

The rest of the paper is organized as follows. 
The Born and regularized inverse Born approximation are investigated in Section \ref{sec:K1} while the numerical results indicating the capability of the proposed method to recover sufficiently small perturbations $\eta$ and $\gamma$ are presented in Section \ref{sec:num_K1}.
The Born and regularized inverse Born series are constructed in Section \ref{sec:IBS} and the convergence analysis of the IBS is established in Section \ref{sec:conv}.
Specifically, sufficient conditions on the convergence of the IBS are given in Theorems \ref{thm:IBS_conv} and \ref{thm:IBS_conv_2} (and Corollary \ref{cor:IBS_conv}), while the approximation error of the IBS with the proposed regularization technique is given in Theorem \ref{thm:IBS_approx_err}.
Preliminary numerical examples of the performance of the IBS are provided in Section \ref{sec:num}.


\section{Born and inverse Born approximation}\label{sec:K1}
\subsection{Born approximation}\label{sec:Born}
In this section, we study the regularized pseudo-inverse of the forward operator for the linearized problem \eqref{eqn:linear_us} which states
\begin{eqnarray}
&&u^s(x;\hat{\theta};k)\approx  
\int_{\Omega}i k \nabla_x G^k(x,y)\cdot \hat{\theta}\,  \gamma(y)e^{ik y\cdot \hat{\theta}} {\rm d}y+\int_{\Omega} G^k(x,y)k^2 \eta(y) e^{ik y\cdot \hat{\theta}} {\rm d}y.\nonumber
\end{eqnarray}
The asymptotic expressions 
\begin{equation}\label{eqn:asym}
\begin{array}{rl}
G^k(x,y)\;=&\frac{e^{i\frac{\pi}{4}}}{\sqrt{8\pi k}}\frac{e^{ik |x|}}{\sqrt{|x|}}e^{-i k \hat{x}\cdot y}+{O}\left(\frac{1}{|x|^{3/2}}\right)\qquad \qquad \mbox{and}\\
\nabla_x G^k(x,y)\;=&-\nabla_y G^k(x,y)=\frac{e^{i\frac{\pi}{4}}}{\sqrt{8\pi k}}\frac{e^{ik |x|}}{\sqrt{|x|}} ik \hat x e^{-i k \hat{x}\cdot y}+{O}\left(\frac{1}{|x|^{3/2}}\right)
\end{array}
\end{equation}
together with the definition of the far-field pattern $u^\infty(\hat{x},\hat{\theta};k)$ (i.e. the scattering data), imply that
\begin{align*}
u^\infty(\hat{x},\hat{\theta};k)\approx \frac{e^{i\frac{\pi}{4}}k^2}{\sqrt{8\pi k}}  \left(\int_{\Omega} - (\hat{x}\cdot \hat{\theta})   e^{ik (\hat{\theta}-\hat{x})\cdot y}\gamma(y){\rm d}y+\int_{\Omega}  e^{ik (\hat{\theta}-\hat{x})\cdot y}\eta(y){\rm d}y\right).
\end{align*}
The support $\Omega$ of the perturbations is part of the unknown. 
Without loss of generality, we assume that the perturbations are supported in the unit disk centered at the origin $B:=B(0,1)\subset\mathbb{R}^2$, i.e., $\overline \Omega\subset B$. 
Note by the scaling properties of the Helmholtz equation, the given formulation in $\overline \Omega\subset B(0,R)$, for any $R>0$, at frequency $k$ is the same the problem formulated in $\overline \Omega\subset B$ with $k/R$ being the frequency.
Then by extending $\gamma$ and $\eta$ by zero in $B \setminus \overline{\Omega}$, we arrive at the following expression for the approximate scattering data
\begin{align*}
u_b^\infty(\hat{x},\hat{\theta};k):=
&\frac{e^{i\frac{\pi}{4}}k^2}{\sqrt{8\pi k}} \left(\int_{B} - (\hat{x}\cdot \hat{\theta})   e^{ik (\hat{\theta}-\hat{x})\cdot y}\gamma(y){\rm d}y+\int_{B}  e^{ik (\hat{\theta}-\hat{x})\cdot y}\eta(y){\rm d}y\right).
\end{align*}
Obviously, the extension by zero of $\gamma$ in $B$ is still piecewise continuously differentiable.

It follows from a direct calculation from the fact $\hat x, \hat{\theta}\in{\mathbb S}$ that $
\hat{\theta}-\hat{x} \in 2\overline{B}$. 
We introduce a new variable
$p := \frac{1}{2}(\hat{\theta}-\hat{x})\in \overline{B}$ and consider the scaled approximate scattering data as a function of $p$, defined by
\begin{align*}
u_b(p;k):=&\sqrt{8\pi}e^{-i\frac{\pi}{4}}k^{-\frac{3}{2}}u_b^\infty(\hat{x},\hat{\theta};k)=\int_{B} - (\hat{x}\cdot \hat{\theta})   e^{ik (\hat{\theta}-\hat{x})\cdot y}\gamma(y){\rm d}y+\int_{B}  e^{ik (\hat{\theta}-\hat{x})\cdot y}\eta(y){\rm d}y\\
=& (2|p|^2-1)\int_{B} e^{i2kp\cdot y}\gamma(y) {\rm d}y+\int_{B} e^{i2kp\cdot y}\eta(y){\rm d}y,
\end{align*}
where the last equation is derived using the identity
$$-(\hat{x}\cdot \hat{\theta}) = \frac{|\hat{\theta}- \hat{x}|^2 - |\hat{\theta}|^2 - |\hat{x}|^2}{2} = 2|p|^2-1.$$
By defining the restricted Fourier operator
\begin{equation}\label{eqn:F} 
\mathcal{F}^k: \;\;\;\;f \mapsto \mathcal{F}^k(p;f) =\int_{B} e^{i2kp\cdot y}  f(y) \ind y, \quad  p\in B,
\end{equation}
we obtain 
\begin{equation}\label{eqn:data u_b(p;k) formula}
u_b(p;k)=(2|p|^2-1)\mathcal{F}^k(p;\gamma) +\mathcal{F}^k(p;\eta), \quad p\in B.
\end{equation}

To reconstruct the perturbations $\gamma$ and $\eta$, we need to solve the linearized ill-posed inverse problem \eqref{eqn:linear_us} with scaled scattering data at two different wave frequencies $k$ and $\ell k$ with $\ell>1$, i.e., $\{u_b(p;k) \cup u_b(p;\ell k)\;| \; p \in B\}$
where
\begin{align*}
u_b(p;k) =& (2|p|^2-1)\mathcal{F}^k (p;\gamma)+ \mathcal{F}^k (p;\eta),\\
\mbox{and}\quad u_b(p;\ell k) =& (2|p|^2-1)\mathcal{F}^{\ell k} (p;\gamma)+ \mathcal{F}^{\ell k} (p;\eta)\\
=&(2|p|^2-1) \int_B e^{i 2\ell k p\cdot y}\gamma(y) \ind y + \int_B e^{i 2\ell k p\cdot y}\eta(y)  \ind y\\
=&(2|p|^2-1)\mathcal{F}^k (\ell p;\gamma)+ \mathcal{F}^k(\ell p;\eta), 
\end{align*}
which yields
\begin{equation}\label{eqn:data u_b(ell k) formula}
    u_b(\ell^{-1}p; \ell k) = (2\ell^{-2}|p|^2-1)\mathcal{F}^k (p;\gamma)+ \mathcal{F}^k (p;\eta) , \quad p \in B.
\end{equation}
Combining formulas \eqref{eqn:data u_b(p;k) formula} and \eqref{eqn:data u_b(ell k) formula}, we derive the following formula of the scattering data at two wave frequencies $k$ and $\ell k$ in the matrix form, i.e.,
\begin{equation}\label{eqn:K1}
\begin{bmatrix}
u_b(p;k)\\
u_b(\ell^{-1}p; \ell k) 
\end{bmatrix}=\begin{bmatrix}
2|p|^2-1&1\\
2\ell^{-2}|p|^2-1&1
\end{bmatrix}\mathcal{F}^k\bigg(p;\begin{bmatrix}
\gamma\\
\eta
\end{bmatrix}\bigg):=A(p)\mathcal{F}^k\bigg(p;\begin{bmatrix}
\gamma\\
\eta
\end{bmatrix}\bigg):=K_b\bigg(\begin{bmatrix}
\gamma\\
\eta
\end{bmatrix}\bigg),
\end{equation} 
where $\mathcal{F}^k$ is applied element-wise to the vector-valued functions.

\subsection{Regularized inverse Born approximation}\label{sec: disk PSWFs}
When the operator $\mathcal{F}^k$ is invertible and $0\neq p \in B$, the above matrix equation \eqref{eqn:K1} lead to the following reconstruction formula
\begin{align*}
\begin{bmatrix}
\gamma\\
\eta
\end{bmatrix}=K_b^{-1}\bigg(\begin{bmatrix}
u_b(p;k)\\
u_b(\ell^{-1}p; \ell k) 
\end{bmatrix}\bigg)=(\mathcal{F}^k)^{-1}\bigg(A(p)^{-1}\begin{bmatrix}
u_b(p;k)\\
u_b(\ell^{-1}p; \ell k) 
\end{bmatrix}\bigg),
\end{align*}
where $A(p)^{-1}=\frac{\ell^2}{2|p|^2(\ell^2-1)}\begin{bmatrix}
1&-1\\
1-2\ell^{-2}|p|^2&2|p|^2-1
\end{bmatrix}$ and $(\mathcal{F}^k)^{-1}$ is applied element-wise to the vector-valued functions.
We consider the regularized pseudo-inverses of $\mathcal{F}^k$ and $A(p)$, denoted by $(\mathcal{F}^k)^\dag$ and $A^\dag(p)$ respectively, when $\mathcal{F}^k$ and $A(p)$ are not invertible or the inverses are unbound. 
Then we define the regularized pseudo-inverses of $K_b$ by \begin{equation}\label{eqn:K1_inverse}
K_b^\dag: \;\;\;\;{\bf f} \mapsto K_b^\dag({\bf f}) =(\mathcal{F}^k)^\dag\big(A^\dag(p) {\bf f}(p)\big)
\end{equation}
with ${\bf f}$ being the vector-valued function.
For the inner component $A^\dag(p)$, we set \begin{equation}\label{eqn:A}
A^\dag(p):=\frac{\ell^2}{2\max(\epsilon,|p|)^2(\ell^2-1)}\begin{bmatrix}
1&-1\\
1-2\ell^{-2}|p|^2&2|p|^2-1
\end{bmatrix}    
\end{equation}
for some regularization parament $\epsilon>0$.
For the outer component $\mathcal{F}^k$, it is worth noting that the operator $\mathcal{F}^k$ exhibits a low-rank structure, which motivates the study of $(\mathcal{F}^k)^\dag$.
We refer the reader to \cite{meng23data,ZhouAudibertMengZhang:2024} for its low-rank approximation based on disk prolate spheroidal wave functions (PSWFs), employed in this work, and to \cite{Zhou2coef:2025} for its complementary low-rank approximation using butterfly factorization and neural networks. 
For the reader's convenience, we introduce the inverse Born reconstruction via a spectral cutoff regularization based on the disk PSWFs \cite{meng23data,ZhouAudibertMengZhang:2024} below.

First, we introduce the spectral decomposition of $\mathcal{F}^k$ considered in this work. 
According to \cite{Slepian64}, there exist real-valued  eigenfunctions $\{\psi_{m,n,l}(x;c)\}^{l\in\mathbb{I}(m)}_{m,n\in \mathbb{N}}$ of the restricted Fourier operator $\mathcal{F}^k$, defined in \eqref{eqn:F}, with bandwidth parameter $c:=2k$:
\begin{align}\label{eqn:eigen_R_Fourier}
    \mathcal{F}^k \psi_{m,n,l}(x)
    &=\int_{B}e^{ic x\cdot y}\psi_{m,n,l}(y;c)dy \nonumber\\
    &=\alpha_{m,n}(c)\psi_{m,n,l}(x;c),\quad x\in B
\end{align}
with $\psi_{m,n,l}(x;c)$ being the disk PSWF and $\alpha_{m,n}(c)$ being the corresponding prolate eigenvalue, 
where $\mathbb{N}=\{0,1,2,3,\dots\}$ and 
\begin{eqnarray*}
    \mathbb{I}(m)=\left\{
    \begin{array}{cc}
        \{1\} & m=0 \\
        \{1,2\} & m \geq 1
    \end{array}\right..
\end{eqnarray*} 
Note that the disk PSWFs are also eigenfunctions of a Sturm-Liouville operator \cite{Slepian64}, i.e.,
\begin{equation}\label{sturm-liouvill}
    \mathcal{D}_c [\psi_{m,n,l}](x)=\chi_{m,n} \psi_{m,n,l}(x),\quad x\in B,
\end{equation}
where $\chi_{m,n}(c)$ is the Sturm-Liouville eigenvalue and 
\begin{align*}
    \mathcal{D}_{c} := -(1-r^2)\partial_r^2-\frac{1}{r}\partial_r+3r\partial_r-\frac{1}{r^2}\Delta_0+c^2r^2
\end{align*}
with the Laplace–Beltrami operator $\Delta_0= \partial^2_\theta $ being the spherical part of Laplacian $\Delta$.
Thanks to the above eigensystems, the following lemma \cite{Slepian64,ZLWZ20} indicates a low-rank structure of $\mathcal{F}^k$ and guarantees that the eigenfunctions, i.e., disk PSWFs, can be computed using the Sturm–Liouville operator, ensuring both stability and efficiency. 
We refer the reader to \cite{ZLWZ20} for the explicit algorithm used in the computation.
\begin{lemma}
For any $c>0$, $\{\psi_{m,n,l}(x;c)\}^{l\in\mathbb{I}(m)}_{m,n\in \mathbb{N}}$ forms a complete and orthonormal system of $L^2(B)$, i.e., for $\forall~m,~n,~m',~n'\in\mathbb{N},~l\in\mathbb{I}(m),~l'\in\mathbb{I}(m')$, there holds
\begin{equation*}
    \int_{B}\psi_{m,n,l}(y;c)\psi_{m',n',l'}(y;c)\ind y=\delta_{m m'}\delta_{n n'}\delta_{l l'},
\end{equation*}
where $\delta$ denotes the Kronecker delta.
\begin{enumerate}
\item[{\rm (i)}] The corresponding prolate eigenvalues $\{\alpha_{m,n}(c)\}_{m,n\in \mathbb{N}}$ in \eqref{eqn:eigen_R_Fourier} are non-zero, and $\lambda_{m,n}(c) := |\alpha_{m,n}(c)|$ can be ordered for fixed $m$ as 
\begin{equation*}
    \lambda_{m,n_1}(c)>\lambda_{m,n_2}(c)>0, \quad \forall n_1<n_2.
\end{equation*}
Moreover, $\lambda_{m,n}(c)\longrightarrow 0$ as $m,n\longrightarrow +\infty$.
\item[{\rm (ii)}] The corresponding Sturm-Liouville eigenvalues $\{\chi_{m,n}\}_{m,n\in \mathbb{N}}$ in \eqref{sturm-liouvill} are real positive and are ordered for fixed $m$ as 
\begin{equation*}
    0<\chi_{m,n_1}(c)<\chi_{m,n_2}(c), \quad \forall n_1<n_2.
\end{equation*}
\end{enumerate}
\end{lemma}
\noindent It is worth noting that the prolate eigenvalues  $\alpha_{m,n}(c)$ decay to zero exponentially fast and the dominant prolate eigenvalues are numerically the same.

Now, we can state the spectral cutoff regularization (e.g., \cite{meng23data,ZhouAudibertMengZhang:2024}) for the outer component $(\mathcal{F}^k)^\dag$ in the inverse Born reconstruction based on the above observation. 
In particular, for any $f\in L^2(B)$, we have the following spectral decomposition of $\mathcal{F}^k$:
\begin{equation}\label{eqn:spec_decom}
\mathcal{F}^k (f)= \sum_{m,n,l} \alpha_{m,n}(c)  \left\langle {f},  \psi_{m,n,l}(\cdot;c) \right\rangle_{B}  \psi_{m,n,l}(\cdot;c),
\end{equation}
with $\langle \cdot,\cdot\rangle_{B}$ being the $L^2(B)$ inner product, and define 
\begin{equation}\label{eqn:low-rank IB reconstruction}
 (\mathcal{F}^k)^{\dag} (f)= \sum_{|\alpha_{m,n}(c)|\geq\alpha}   \frac{1}{ \alpha_{m,n}(c) } \left\langle {f},  \psi_{m,n,l}(\cdot;c) \right\rangle_{B}  \psi_{m,n,l}(\cdot;c),
\end{equation}
with $\alpha>0$ being the spectral cutoff parameter.


\section{Born series and inverse Born series}\label{sec:IBS}
In this section, we shall first derive the Born series of the scattered wave field $u^s=\sum_{j=1}^\infty u_j^s$ by using the recursive relation of its sequential terms derived from the Born approximation \eqref{eqn:linear_us}, i.e.,
\begin{align}
u^s_{j+1}(x;\hat{\theta};k)= &\nabla_x\cdot \int_{B} G^k(x,y)\gamma(y)\nabla_y u^s_j(y;\hat{\theta};k) {\rm d}y+\int_{B} G^k(x,y)k^2 \eta(y) u^s_j(y;\hat{\theta};k){\rm d}y\nonumber\\
= &\int_{B}\nabla_x\cdot G^k(x,y)\, \gamma(y) \nabla_y u^s_j(y;\hat{\theta};k) {\rm d}y+k^2\int_{B} G^k(x,y) \eta(y)u^s_j(y;\hat{\theta};k){\rm d}y\nonumber\\
= &\int_{B}\nabla_x G^k(x,y)\gamma(y) \cdot \nabla_y u^s_j(y;\hat{\theta};k) {\rm d}y+k^2\int_{B} G^k(x,y) \eta(y)u^s_j(y;\hat{\theta};k){\rm d}y\nonumber,
\end{align}
with $u^s_0(x;\hat{\theta};k)=u^i(x;\hat{\theta};k)=e^{ik\hat{\theta}\cdot x}$ being the incident wave field.
After defining the integral operators
\begin{eqnarray}
F_0^k:& f \mapsto F_0^k(x;f)=k^2\int_{B} G^k(x,y)f(y) {\rm d}y, \quad  &x\in B\label{eqn:F0}\\
\mbox{and}\quad F_1^k:& \;\;\;\; {\bf f} \mapsto F_1^k(x;{\bf f})=\int_{B}\nabla_{x} G^k(x,y)\cdot{\bf f}(y)\; {\rm d}y, \quad  &x\in B \label{eqn:F1}    
\end{eqnarray}
with $f$ being the scalar-valued functions and ${\bf f}$ being the vector-valued function, we derive the recursion
\begin{align}
u^s_{j+1}(x;\hat{\theta};k)=
&F_1^k(x;\gamma(\cdot)\nabla u^s_j(\cdot;\hat{\theta};k))+F_0^k(x;\eta(\cdot)u^s_j(\cdot;\hat{\theta};k)).\label{eqn:recursion_forward}
\end{align}

Now, we start with the first-order iteration $u_1^s$. The recursion \eqref{eqn:recursion_forward} with $j=0$ and $u_0^s=u^i$ gives
\begin{align*}
u^s_1(x;\hat{\theta};k)
= &\int_{B}\nabla_{x} G^k(x,y)\gamma(y)\cdot  \nabla_{y} u^i(y;\hat{\theta};k) {\rm d}y+k^2\int_{B} G^k(x,y) \eta(y)u^i(y;\hat{\theta};k){\rm d}y\\
=&F_1^k(x;\gamma(\cdot)\nabla u^i(\cdot;\hat{\theta};k))+F_0^k(x;\eta(\cdot) u^i(\cdot;\hat{\theta};k)).
\end{align*}
Further, the gradients of $u^s_1(x)$ with respect to $x$ is given by
\begin{align*}
\nabla_{x} u^s_1(x;\hat{\theta};k)
= &\int_{B}\nabla_{x}\nabla_{x} G^k(x,y) \gamma(y)\;\nabla_{y} u^i(y;\hat{\theta};k)   {\rm d}y+k^2\int_{B} \nabla_{x} G^k(x,y) \eta(y)u^i(y;\hat{\theta};k) {\rm d}y\\
:=&{\bf F_3^k}(x;\gamma(\cdot)\nabla u^i(\cdot;\hat{\theta};k))+ {\bf F_2^k}(x;\eta(\cdot)u^i(\cdot;\hat{\theta};k)),
\end{align*} 
where 
\begin{eqnarray}
{\bf F_2^k}:&  \;\;\;f \mapsto {\bf F_2^k}(x;f)=k^2\int_{B} \nabla_{x}G^k(x,y)f(y) {\rm d}y, \quad & x\in B\label{eqn:F2}\\
\mbox{and}\quad {\bf F_3^k}:& \;\;\;\;{\bf f} \mapsto {\bf F_3^k}(x;{\bf f})=\int_{B}\nabla_{x}\nabla_{x} G^k(x,y){\bf f}(y)\; {\rm d}y, \quad  &x\in B. \label{eqn:F3}    
\end{eqnarray}

We recall the relation $p = \frac{\hat{\theta}-\hat{x}}{2}$ and define $q(p)=\frac{\sqrt{1-|p|^2}}{|p|}\begin{bmatrix}
    0&-1\\
    1&0
\end{bmatrix}p$ such that 
\begin{align*}
q(p)=(-\sqrt{1-|p|^2}\sin\theta,\sqrt{1-|p|^2}\cos\theta)^t  \quad\mbox{for any}\quad p=(|p|\cos\theta,|p|\sin\theta)^t,
\end{align*}
and thus 
\begin{align*}
&\hat{\theta}=q(p)+p \qquad\;\;\mbox{and}\quad  \hat{x}=q(p)-p, \qquad\;\;\mbox{when}\quad \theta_{\hat{x}}\in(\theta_{\hat{\theta}},\theta_{\hat{\theta}}+\pi),\\
&\hat{\theta}=-(q(p)-p) \quad\mbox{and}\quad  \hat{x}=-(q(p)+p), \quad\mbox{when}\quad \theta_{\hat{x}}\in(\theta_{\hat{\theta}}-\pi,\theta_{\hat{\theta}}).
\end{align*}
Thanks to the reciprocity relation of the far-field pattern for the $j$th iterative scattered wave field $u^\infty_j(\hat{x};\hat{\theta};k)=u^\infty_j(-\hat{\theta};-\hat{x};k)$, we only need to consider the case where $\theta_{\hat{x}}\in(\theta_{\hat{\theta}},\theta_{\hat{\theta}}+\pi)$. 
Above identities yield $$u^i(y;p;k)=e^{ik \hat{\theta}\cdot y}=e^{ik (q+p)\cdot y} \quad\mbox{and} \quad\nabla_y u^i(y;p;k)=ik\hat{\theta} e^{ik\hat{\theta}\cdot y}=ik (q+p)u^i(y;p;k),$$
and thus we rewrite the $j$th iterative scattered wave field and its gradient as functions of $p$ and $k$ (and $y$) as
\begin{align*}
u^s_1(x;p;k)=&F_1^k\big(x; \gamma(y)\nabla_y u^i(y;p;k)\big)+F_0^k(x; \eta(y) u^i(y;p;k)),\\
\nabla u^s_1(x;p;k)=&{\bf F_3^k}\big(x; \gamma(y)\nabla_y u^i(y;p;k)\big)+ {\bf F_2^k}(x; \eta(y) u^i(y;p;k)).
\end{align*}
For the sake of simplifying the analysis, we use the notations $u^i(\cdot)=u^i(\cdot;p;k)$ and $u^s_{j}(\cdot)=u^s_j(\cdot;p;k)$ for any $j\geq1$ and consider the scaled far-field pattern of the $j$th iterative scattered wave field defined as $$u_j(p;k):=\sqrt{8\pi}e^{-i\frac{\pi}{4}}k^{-\frac{3}{2}}u_j^\infty(\hat{x};\hat{\theta};k)$$ with the relation $p = \frac{\hat{\theta}-\hat{x}}{2}$. 
Note that the first-order term is the scaled Born-approximated scattering data discussed in Section \ref{sec:Born}, i.e., $u^\infty_1=u^\infty_b$ and $u_1=u_b$.
Then the asymptotic expressions of $G^k(x,y)$ and $\nabla_x G^k(x,y)$ in \eqref{eqn:asym} imply that
\begin{align*}
u_1(p;k)
=&k^{-2} \int_{B} ik\hat{x} e^{-ik\hat{x}\cdot y}\gamma(y)\cdot \nabla u^i(y) {\rm d}y+\int_{B} e^{-ik\hat{x}\cdot y}\eta(y) u^i(y){\rm d}y\\
=&-\int_{B}  e^{-ik(q-p)\cdot y}\gamma(y)(q-p)\cdot (q+p) u^i(y) {\rm d}y+\int_{B} e^{-ik(q-p)\cdot y}\eta(y) u^i(y){\rm d}y\\
:=&\mathcal{F}^{p,k}_0\big(\gamma (p-q)\cdot (q+p)u^i\big)
+\mathcal{F}^{p,k}_0(\eta u^i),
\end{align*}
where $\mathcal{F}^{p,k}_0(f)=\int_{B}e^{ik(p-q)\cdot y}f(y){\rm d}y$.
Furthermore, we can derive that
\begin{align*}
u_1(p;k)
=&u_b(p;k) 
= (2|p|^2-1)\mathcal{F}^k(p;\gamma)+\mathcal{F}^k(p;\eta),\\ 
u_1(\ell^{-1} p;\ell k) =& u_b(\ell^{-1} p;\ell k)
=(2\ell^{-2}|p|^2-1)\mathcal{F}^k(p;\gamma)+\mathcal{F}^k(p;\eta).
\end{align*}
We define the operator $$K_1(\cdot):=K_b(\cdot)=\begin{bmatrix}
2|p|^2-1&1\\
2\ell^{-2}|p|^2-1&1
\end{bmatrix}\mathcal{F}^k(p;\cdot),$$ 
where $\mathcal{F}^k$ defined in \eqref{eqn:F} is applied to the vector-valued function element-wise, such that $K_1\big((\gamma,\eta)^t\big)=\big(u_1(p;k),u_1(\ell^{-1}p;\ell k)\big)^t$.

Next, we consider the second iteration $u_2^s$. 
Together with the above estimates on $u^s_1$ and its gradient, we derive from the recursion \eqref{eqn:recursion_forward} and the relations $\nabla F_1^k={\bf F_3^k}$ and $\nabla F_0^k={\bf F_2^k}$ that
\begin{align*}
u^s_{2}
=& F_1^k(\gamma\nabla u_1^s)+F_0^k(\eta u_1^s)\\
=&F_1^k\Big(\gamma \big({\bf F_3^k}(\gamma\nabla u^i)+{\bf F_2^k}(\eta u^i)\big)\Big)+F_0^k\Big(\eta \big(F_1^k(\gamma\nabla u^i)+F_0^k(\eta  u^i)\big)\Big),\\
\nabla u^s_{2}
=&{\bf F_3^k}(\gamma\nabla u^s_{1})+{\bf F_2^k}(\eta u^s_{1}),
\end{align*}
and thus
\begin{align*}
u_2(p;k)
= &k^{-2} \int_{B} ik\hat{x} e^{-ik\hat{x}\cdot y}\gamma(y)\cdot \nabla u_1^s(y) {\rm d}y+\int_{B} e^{-ik\hat{x}\cdot y}\eta(y) u_1^s(y){\rm d}y\\
= &ik^{-1} \int_{B} (q-p) e^{ik(p-q)\cdot y}\gamma(y)\cdot \nabla u_1^s(y) {\rm d}y+\int_{B} e^{ik(p-q)\cdot y}\eta(y) u_1^s(y){\rm d}y\\
= &ik^{-1} \mathcal{F}^{p,k}_0(\gamma (q-p)\cdot \nabla u_1^s)+\mathcal{F}^{p,k}_0(\eta u_1^s)\\
= &ik^{-1} \mathcal{F}^{p,k}_0\Big(\gamma (q-p)\cdot \big({\bf F_3^k}(\gamma,\nabla u^i)+{\bf F_2^k}(\eta,u^i)\big)\Big)+\mathcal{F}^{p,k}_0\Big(\eta \big(F_1^k(\gamma\nabla u^i)+F_0^k(\eta  u^i)\big)\Big).
\end{align*}
Then, we define the operator $K_2$ by $K_2\big((\gamma,\eta)^t, (\gamma,\eta)^t\big)=\big(u_2(p;k),u_2(\ell^{-1}p;\ell k)\big)^t$.

Similarly, for any $j\geq 1$, we can derive that
\begin{align}
u^s_{j}=&F_1^k(\gamma\nabla u^s_{j-1})+F_0^k(\eta u^s_{j-1}),\quad
\nabla u^s_{j}={\bf F_3^k}(\gamma \nabla u^s_{j-1})+{\bf F_2^k}(\eta u^s_{j-1}),\label{eqn:us_j}\\
u_j(p;k)=&ik^{-1}\mathcal{F}^{k,p}_0(\gamma (q-p)\cdot\nabla u^s_{j-1})+\mathcal{F}^{k,p}_0(\eta u^s_{j-1}),\label{eqn:u_j}
\end{align}
and we define the $j$-multilinear operator $K_j$ such that $$K_j\big((\gamma,\eta)^t, \cdots, (\gamma,\eta)^t\big)=\big(u_j(p;k),u_j(\ell^{-1}p;\ell k)\big)^t.$$ 

Let $\phi=\sum_{j=1}^{\infty}(u_j(p;k),u_j(\ell^{-1}p;\ell k)\big)^t$ and $\psi=(\gamma,\eta)^t$, then the Born series for the scaled far-field pattern $\phi$ is given by
\begin{align*}
\phi= \sum_{j=1}^{\infty} K_j\big(\psi, \cdots,  \psi\big).
\end{align*}
Given the above Born series for $\phi$, the IBS \cite{HoskinsSchotland:2022} is defined as
\begin{align}\label{eqn:IBS}
\psi:=\sum_{j=1}^{\infty} \psi_j:=\sum_{j=1}^{\infty} \mathcal{K}_j(\phi),
\end{align}
with
\begin{align*}
\psi_1:=&\mathcal{K}_1(\phi)=K_1^\dag(\phi),\\
\psi_2:=&\mathcal{K}_2(\phi)=-\mathcal{K}_1\big(K_2(\mathcal{K}_1(\phi),\mathcal{K}_1(\phi))\big)
=-\mathcal{K}_1\big(K_2(\psi_1,\psi_1)\big),\\
\psi_3:=&\mathcal{K}_3(\phi)=-\mathcal{K}_1\big(K_2(\psi_1, \psi_2)+K_2(\psi_2, \psi_1)+K_3(\psi_1, \psi_1,\psi_1)\big),\\
\vdots &\\
\psi_j:=&\mathcal{K}_j(\phi)=-\mathcal{K}_1\big(\sum_{m=2}^{j}\sum_{\sum_{t=1}^m i_t=j}K_m(\psi_{i_1},\cdots, \psi_{i_m})\big),
\end{align*}
where $\mathcal{K}_1=K_1^\dag=K_b^\dag$, defined in \eqref{eqn:K1_inverse}, is the regularized pseudo-inverse of the operator $K_1$.


\section{Convergence of the regularized inverse Born series}\label{sec:conv}
In this section, we discuss the sufficient conditions for convergence of the IBS \eqref{eqn:IBS} and its approximation error. 

\subsection{General regularized inverse Born series}\label{conv_gen}
We start with the convergence of the general regularized inverse Born series with a regularized and bounded operator $\mathcal{K}_1$.
First, we state a sufficient condition for the convergence of the general IBS established in \cite[Theorem 2.2]{HoskinsSchotland:2022}. 
We also refer the reader to the stability of the general IBS in \cite[Theorem 3.2]{MoskowSchotland:2008}, where the estimates depend on $\nu$, $\mu$ and $\|\mathcal{K}_1\|_{(L^2(B))^2\to \mathcal{X}^2}$ introduced in the following lemma.
\begin{lemma}\label{lem:conv_IBS_Hoskins}
Let $\mu$ and $\nu$ be positive constants such that 
\begin{equation*}
\|K_j(\tilde{\psi}_1,\cdots,\tilde{\psi}_j)\|_{(L^2(B))^2}\leq  \nu \mu^{j-1}\Pi_{i=1}^j \|\tilde{\psi}_j\|_{\mathcal{X}^2}, \quad \mbox{for any} \quad j\geq 1,
\end{equation*}
with $\mathcal{X}$ being a Banach space and $\|\textbf{v}\|_{\mathcal{X}^2}:=(\sum_{j=1}^2 \|v_j\|^2_{\mathcal{X}})^\frac12
$, for any vector-valued function $\textbf{v}=(v_1,v_2)^t 
\in {\mathcal{X}}^2$. 
The IBS \eqref{eqn:IBS} converges if $\|\mathcal{K}_1 (\phi)\|_{{\mathcal{X}}^2}<r$ with the radius of convergence $r=
\big(2\mu(\sqrt{16C^2+1}+4C)\big)^{-1}$, where $C=\max(2,\nu\|\mathcal{K}_1\|_{(L^2(B))^2\to {\mathcal{X}}^2})$.
\end{lemma}

Now, we shall specify the constants $\nu$ and $\mu$ in the following proposition and defer its proof to Appendix \ref{app:estimate} and the estimation of the corresponding coefficients $a(k)$, $b(k)$ and $c(k)$ to Lemma \ref{lem:abc}.
\begin{proposition}\label{prop:nu_mu}
Let $\tilde{\phi}_j:=K_j(\tilde{\psi}_1,\cdots,\tilde{\psi}_j)$, then $\|\tilde{\phi}_{j}\|_{(L^2(B))^2}\leq \nu_{\infty} \mu_{\infty}^{j-1}\Pi_{i=1}^{j} \|\tilde{\psi}_i\|_{(L^\infty(B))^2}$ with
$\nu_{\infty}=\sqrt{2}|B|$ and $\mu_{\infty} = \sqrt{2} \big(\mu_0(k)+\mu_0(\ell k)\big)$, where $$\mu_0(k)=\max(1,k b(k), k^{-1}c(k),k^2 |B|^\frac12 a(k)),$$ $a(k)=\sup_{x\in B}\|G^k(x,\cdot)\|_{L^2(B)}$, $b(k)=\|F_1^k\|_{(L^2(B))^2\to L^2(B)}$ and $c(k)=\|{\bf F_2^k}\|_{L^2(B)\to (L^2(B))^2}$.
\end{proposition}

\begin{corollary}\label{cor:bound_L2}
Let $\tilde{\phi}_j:=K_j(\tilde{\psi}_1,\cdots,\tilde{\psi}_j)$ and \begin{equation}\label{eqn:measure}
M=\min\Big(\big|\{x\in B: \;\gamma(x)\geq\frac12\|\gamma\|_{L^\infty(B)}\}\big|,\big|\{x\in B: \;\eta(x)\geq\frac12\|\eta\|_{L^\infty(B)}\}\big|\Big),    
\end{equation} 
then $\|\tilde{\phi}_{j}\|_{(L^2(B))^2}\leq \nu_{2} \mu_{2}^{j-1}\Pi_{i=1}^{j} \|\tilde{\psi}_i\|_{(L^2(B))^2}$ with
$\nu_{2}=2 M^{-\frac12}\nu_{\infty}$ and $\mu_{2}=2 M^{-\frac12}\mu_{\infty}$.
\end{corollary}
\begin{proof}
It is well-known that $M\in(0,|B|)$. Then it follows directly from Proposition \ref{prop:nu_mu} and the estimate
\begin{align*}
\|\tilde{\psi}_i\|_{(L^2(B))^2}\geq \frac12 M^\frac12 \|\tilde{\psi}_i\|_{(L^\infty(B))^2},
\end{align*}
that
\begin{align*}
\|\tilde{\psi}_i\|_{(L^\infty(B))^2}\leq 2 M^{-\frac12}\|\tilde{\psi}_i\|_{(L^2(B))^2}.
\end{align*}
\end{proof}

Now, we can give the convergence conditions of the IBS \eqref{eqn:IBS} from Lemma \ref{lem:conv_IBS_Hoskins} with the help of Proposition \ref{prop:nu_mu} and Lemma \ref{lem:abc}.
\begin{theorem}\label{thm:IBS_conv}
Let $k>\frac12$ and $\|\mathcal{K}_1\|_{(L^2(B))^2\to (L^\infty(B))^2}\leq \tau^{-1}$ with $\tau>0$, then the IBS \eqref{eqn:IBS} converges if $\|\mathcal{K}_1 (\phi)\|_{(L^\infty(B))^2}<r$ with the radius of convergence $$r
\geq c_{r,\infty}(1+\ell^\frac32)^{-1}k^{-\frac32}\min(2^{-\frac12}\pi,\tau)$$ with a constant $c_{r,\infty}$ that is independent of $\ell$, $k$ or $\tau$. 
\end{theorem}
\begin{proof}
It follows directly from Lemma \ref{lem:abc}, the definitions of $\mu_0(k)$ and $\mu_{\infty}$ in Proposition \ref{prop:nu_mu} and the assumptions $|B|=\pi$ and $k>\frac12$ that 
\begin{align*}
\mu_0(k)=&\max(1,k b(k), k^{-1}c(k),k^2 |B|^\frac12 a(k))
\leq \big(3\sqrt{2\pi}+2\sqrt{\frac{2}{3}}\pi\big) k^\frac32,\\
\mu_{\infty}=&\sqrt{2} \big(\mu_0(k)+\mu_0(\ell k)\big)
\leq \big(6\sqrt{\pi}+\frac{4\pi}{\sqrt{3}}\big)  (1+\ell^\frac32) k^\frac32.
\end{align*}
Thus, by Lemma \ref{lem:conv_IBS_Hoskins} and the assumption $\|\mathcal{K}_1\|_{(L^2(B))^2\to (L^\infty(B))^2}\leq\tau^{-1}$, the radius of convergence $r$ satisfies
\begin{align*}
r=&\big(2\mu_{\infty}(\sqrt{16C_\infty^2+1}+4C_\infty)\big)^{-1}
\geq \Big( \big(6\sqrt{\pi}+\frac{4\pi}{\sqrt{3}}\big) (1+\ell^\frac32) k^\frac32(8C_\infty+1)\Big)^{-1}\\
\geq &\Big(18 \big(3\sqrt{\pi}+\frac{2\pi}{\sqrt{3}}\big) (1+\ell^\frac32) k^\frac32 C_\infty\Big)^{-1}
\end{align*} 
where 
\begin{align*}
C_\infty=&\max(2,\nu_{\infty}\|\mathcal{K}_1\|_{(L^2(B))^2\to (L^\infty(B))^2})\\
=&\max(2,\sqrt{2}|B|\|\mathcal{K}_1\|_{(L^2(B))^2\to (L^\infty(B))^2})
\leq \sqrt{2}\pi\min(2^{-\frac12}\pi,\tau)^{-1},
\end{align*}
which imply
\begin{align*}
r 
\geq& \Big(18\sqrt{2} \big(3\sqrt{\pi}+\frac{2\pi}{\sqrt{3}}\big) \pi  \Big)^{-1}(1+\ell^\frac32)^{-1}k^{-\frac32}\min(2^{-\frac12}\pi,\tau).
\end{align*} 
This completes the proof.
\end{proof}

\begin{corollary}\label{cor:IBS_conv}
Let $k>\frac12$, $M$ be defined in \eqref{eqn:measure}, and $\|\mathcal{K}_1\|_{(L^2(B))^2\to (L^2(B))^2}\leq \tau^{-1}$ with $\tau>0$, then the IBS \eqref{eqn:IBS} converges if $\|\mathcal{K}_1 (\phi)\|_{(L^2(B))^2}<r$ with the radius of convergence $$r
\geq c_{r,2} (1+\ell^\frac32)^{-1}k^{-\frac32} M^\frac12\min(\sqrt{2}\pi ,\tau M^\frac12)$$ with a constant $c_{r,2}$ that is independent of $\ell$, $k$, $M$ or $\tau$.  
\end{corollary}
\begin{proof}
Following the analysis in the proof of Theorem \ref{thm:IBS_conv}, the relations $\nu_{2}=2 M^{-\frac12}\nu_{\infty}$ and $\mu_{2}=2 M^{-\frac12}\mu_{\infty}$ in Corollary \ref{cor:bound_L2} and the estimates in Lemmas \ref{lem:abc} and \ref{lem:conv_IBS_Hoskins} yield that 
\begin{align*}
C_2=&\max(2,\nu_{2}\|\mathcal{K}_1\|_{(L^2(B))^2\to (L^2(B))^2})\leq 2^\frac32 \pi \min(\sqrt{2}\pi ,\tau M^\frac12)^{-1}.
\end{align*} 
Together with the fact that $2 M^{-\frac12}>1$, there holds
\begin{align*}
r \geq& \Big(72 \sqrt{2}\big(3\sqrt{\pi}+\frac{2\pi}{\sqrt{3}}\big) \pi  \Big)^{-1} (1+\ell^\frac32)^{-1}k^{-\frac32} M^\frac12\min(\sqrt{2}\pi ,\tau M^\frac12),
\end{align*} 
which completes the proof.
\end{proof}

\subsection{Regularized inverse Born series with the disk PSWFs}\label{conv_PSWF}
In this section, we study the convergence and approximation error of the IBS \eqref{eqn:IBS} with the proposed regularization technique discussed in Section \ref{sec: disk PSWFs}. 
Numerically, we employ the spectral cutoff regularization for the outer component $(\mathcal{F}^k)^\dag$ of the inverse Born reconstruction $\mathcal{K}_1$, where the input scattering data will be projected onto the disk PSWFs whose corresponding prolate eigenvalues are larger than a spectral cutoff parameter $\alpha=\tilde\alpha\|\mathcal{F}^k\|_{L^2(B)\to L^2(B)}$ with the constant $\tilde\alpha\in(0,1)$, as defined in \eqref{eqn:low-rank IB reconstruction}. 
In addition, let the inner component $A^\dag(p)$ be defined in \eqref{eqn:A}. 
This yields an upper bound of $\|\mathcal{K}_1\|_{(L^2(B))^2\to (L^2(B))^2}$ and the corresponding radius of convergence for IBS \eqref{eqn:IBS} is given in the following theorem.
\begin{theorem}\label{thm:IBS_conv_2}
Let $k>\frac12$, $M$ be defined in \eqref{eqn:measure}, $(\mathcal{F}^k)^\dag$ be defined in \eqref{eqn:low-rank IB reconstruction} with $c=2k$ and $\alpha=\tilde\alpha\|\mathcal{F}^k\|_{L^2(B)\to L^2(B)}$ being the cutoff parameter, and $A^\dag(p)$ be defined in \eqref{eqn:A} for some constants $\tilde{\alpha}\in(0,1)$ and $\epsilon>0$. 
Then there holds
\begin{align*}
\|\mathcal{K}_1\|_{(L^2(B))^2\to(L^2(B))^2}
\leq \frac{k}{2\tilde{\tau}\sqrt{\min(k,2)}}, \quad\mbox{where}\quad \tilde{\tau}=\tilde\alpha \epsilon^2(1-\ell^{-2}).
\end{align*}
Moreover, the radius of convergence $r$ of the IBS \eqref{eqn:IBS} satisfies $$r
\geq c_{r,2}(1+\ell^\frac32)^{-1}k^{-\frac32} M^\frac12\min(\sqrt{2}\pi,2\sqrt{\min(k,2)} k^{-1} \tilde{\tau}M^\frac12)$$ with a constant $c_{r,2}$ that is independent of $\ell$, $k$, $M$ or $\tilde\tau$. 
\end{theorem}

\begin{proof}
By the definition of $\mathcal{K}_1=K_b^\dag$ in \eqref{eqn:K1_inverse}, for any ${\bf f}\in (L^2(B))^2$, there holds
\begin{align*}
\|\mathcal{K}_1 ({\bf f})\|_{(L^2(B))^2}\leq& \|(\mathcal{F}^k)^\dag\|_{L^2(B)\to L^2(B)}\|A^\dag {\bf f}\|_{(L^2(B))^2}\\
\leq& \|(\mathcal{F}^k)^\dag\|_{L^2(B)\to L^2(B)}\sup_{p\in B}\|A^\dag(p)\|_F \|{\bf f}\|_{(L^2(B))^2}\\
\leq& (\tilde\alpha\|\mathcal{F}^k\|_{L^2(B)\to L^2(B)})^{-1} \frac{\ell^2}{\epsilon^2(\ell^2-1)}\|{\bf f}\|_{(L^2(B))^2},
\end{align*}
which implies
\begin{align*}
\|\mathcal{K}_1\|_{(L^2(B))^2\to(L^2(B))^2}
\leq (\tilde\alpha\|\mathcal{F}^k\|_{L^2(B)\to L^2(B)})^{-1} \frac{\ell^2}{\epsilon^2(\ell^2-1)}:=(\tilde{\tau}\|\mathcal{F}^k\|_{L^2(B)\to L^2(B)})^{-1},
\end{align*}
where $\tilde{\tau}=\tilde\alpha \epsilon^2(1-\ell^{-2})$.
Now, we derive a lower bound for  $\|\mathcal{F}^k\|_{L^2(B)\to L^2(B)}\geq d(k)$ which yields
\begin{align*}
\|\mathcal{K}_1\|_{(L^2(B))^2\to(L^2(B))^2}\leq(\tilde{\tau} d(k))^{-1}:=\tau^{-1}.
\end{align*}
In fact, let $f\equiv |B|^{-\frac12}=\pi^{-\frac12}$ such that $\|f\|_{L^2(B)}=1$, then
\begin{align*}
\|\mathcal{F}^k(p;f)\|_{L^2(B)}=&\pi^{-\frac12}\|\int_{B} e^{i2kp\cdot y} {\rm d}y\|_{L^2(B)}.
\end{align*}
By using the Jacobi–Anger expansion \cite[(3.112)]{colton-kress}, we can further decompose $e^{i2kp\cdot y}$ as
\begin{align}\label{eqn:Jacobi–Anger}
e^{i2kp\cdot y}=J_0(2k|p||y|)+2\sum_{n=1}^\infty i^{n} J_n(2k|p||y|) \cos\big(n (\theta_y-\theta_p)\big),
\end{align}
with $\theta_{y}$ and $\theta_p$ denoting the arguments of $y$ and $p$ respectively.
Then we derive that
\begin{align*}
\int_{B} e^{i2kp\cdot y}{\rm d}y=&\int_{B} J_0(2k|p||y|){\rm d}y =2\pi \int_{0}^1 J_0(2k|p|\rho)\rho{\rm d}\rho\\
=&\frac{\pi}{2} k^{-2} |p|^{-2} \int_{0}^{2k|p|} \rho J_0(\rho){\rm d}\rho
=\pi k^{-1} |p|^{-1} J_1(2k|p|),
\end{align*}
and thus
\begin{align*}
\pi\|\mathcal{F}^k(p;f)\|^2_{L^2(B)}=&\|\int_{B} e^{i2kp\cdot y} {\rm d}y\|_{L^2(B)}^2
=\|\pi k^{-1} |p|^{-1} J_1(2k|p|)\|_{L^2(B)}^2\\
=&\pi^2 k^{-2}\int_{B}|p|^{-2}|J_1(2k|p|)|^2 {\rm d}p
=2\pi^3 k^{-2}\int_{0}^1 \rho^{-1}|J_1(2k\rho)|^2 {\rm d}\rho\\
=&2\pi^3 k^{-2}\int_{0}^{2k} \rho^{-1}|J_1(\rho)|^2 {\rm d}\rho
\geq 2\pi^3 k^{-2} \frac{2}{\pi^2}\min(k,2)
\geq 4\pi k^{-2} \min(k,2)
\end{align*}
which implies 
$
\|\mathcal{F}^k\|_{(L^2(B))^2\to (L^2(B))^2}\geq 2\sqrt{\min(k,2)} k^{-1}:=d(k)
$
and thus $$\tau= \tilde{\tau} d(k)= 2\sqrt{\min(k,2)} k^{-1}\tilde{\tau}.$$
Finally, Corollary \ref{cor:IBS_conv} completes the proof.
\end{proof}

Now, we derive the approximation error of the IBS \eqref{eqn:IBS} based on the following lemma established in \cite[Theorem 2.4]{HoskinsSchotland:2022}.

\begin{lemma}\label{lem:IBS_approx_err_Hoskins}
Let the assumptions in Theorem \ref{thm:IBS_conv_2} hold and that the forward and IBS \eqref{eqn:IBS} converge. 
Let $\tilde{\psi}$ denote the sum of the IBS. We assume that $$\mathcal{M}:=\max(\|\psi\|_{(L^2(B))^2},\|\tilde\psi\|_{(L^2(B))^2})\leq \mu_2^{-1} \big(1-\sqrt{1-(1+C_{\mathcal{K}_1})^{-1}}\big),$$ where 
$C_{\mathcal{K}_1}:=\nu_2\|\mathcal{K}_1\|_{(L^2(B))^2\to(L^2(B))^2}$. 
Then the approximation error can be bounded by
\begin{align*}
\|\psi-\sum_{j=1}^{N} \mathcal{K}_j(\phi)\|_{(L^2(B))^2}\leq &2\mu_2 \big(\sqrt{16C_2^2+1}(1-C_{{\rm ratio}})\big)^{-1}C_{{\rm ratio}}^{N+1}\\
&+\Big(1+\big(1-(1-\mu_2\mathcal{M})^{-2}\big)C_{\mathcal{K}_1}\Big)^{-1}\|(I-\mathcal{K}_1K_1)\psi\|_{(L^2(B))^2}
\end{align*}
where $C_2=\max(2,C_{\mathcal{K}_1})$ and $C_{{\rm ratio}}=\|\mathcal{K}_1\phi\|_{(L^2(B))^2}/r$.
\end{lemma}

To provide a more specific upper bound on the approximation error of the IBS \eqref{eqn:IBS}, we evaluate the ability of the proposed regularized inverse Born approximation to recover $\gamma$ and $\eta$ from the Born scattering data in the following lemma.
\begin{lemma}\label{lem:K1_approx_err}
Let the assumptions in Theorem \ref{thm:IBS_conv_2} hold and $\psi=(\gamma,\eta)^t\in (L^2(B))^2$. Then for any $\gamma$ and $\eta$ such that $\big(\|\gamma-\gamma^\alpha\|^2_{L^2(B)}+\|\eta-\eta^\alpha\|^2_{L^2(B)}\big)^\frac12\leq\delta_\alpha$ for some $\gamma^\alpha,\eta^\alpha\in \mbox{span}\{ \psi_{m,n,l}(\cdot;2k): |\alpha_{m,n}(2k)| \geq \alpha\}$ with $\psi_{m,n,l}$ being the disk PSWF discussed in Section \ref{sec: disk PSWFs}
, there holds
\begin{align*}
\|(I-\mathcal{K}_1K_1)\psi\|_{(L^2(B))^2}\leq 3^{-\frac12}\pi\alpha^{-1} \epsilon \|\psi\|_{(L^2(B))^2}+\delta_\alpha.
\end{align*}
\end{lemma}
\begin{proof}
The definitions of $\mathcal{K}_1=K_b^\dag$ in \eqref{eqn:K1_inverse} and $K_1=K_b$ in \eqref{eqn:K1} yield
\begin{align*}
\mathcal{K}_1K_1\psi=(\mathcal{F}^k)^\dag \big(A^\dag(p)A(p)\mathcal{F}^k(p;\psi)\big)=(\mathcal{F}^k)^\dag\big(\frac{|p|^2}{\max(\epsilon,|p|)^2}\mathcal{F}^k(p;\psi)\big).
\end{align*}
Thus, we can decompose $\|(I-\mathcal{K}_1K_1)\psi\|_{(L^2(B))^2}$ into two components
\begin{align*}
&\|(I-\mathcal{K}_1K_1)\psi\|_{(L^2(B))^2}\\
\leq& \|\big(I-(\mathcal{F}^k)^\dag\mathcal{F}^k\big)\psi\|_{(L^2(B))^2}+\|(\mathcal{F}^k)^\dag\Big(\big(\frac{|p|^2}{\max(\epsilon,|p|)^2}-1\big)\mathcal{F}^k(p;\psi)\Big)\|_{(L^2(B))^2}\\
:=&{\rm II}_1+{\rm II}_2.
\end{align*}
For the first term ${\rm II}_1$, the spectral decomposition of $\mathcal{F}^k$ given in \eqref{eqn:spec_decom} with $c=2k$ and the reconstruction formula \eqref{eqn:low-rank IB reconstruction} imply
\begin{align*}
{\rm II}^2_1=& \|\big(I-(\mathcal{F}^k)^\dag\mathcal{F}^k\big)\gamma\|^2_{L^2(B)}+\|\big(I-(\mathcal{F}^k)^\dag\mathcal{F}^k\big)\eta\|^2_{L^2(B)}\\
\leq&\|\gamma-\gamma^\alpha\|^2_{L^2(B)}+\|\eta-\eta^\alpha\|^2_{L^2(B)}
\leq \delta_\alpha^2.
\end{align*}
For the second term ${\rm II}_2$, by the estimate $\|(\mathcal{F}^k)^\dag\|_{L^2(B)\to L^2(B)}\leq \alpha^{-1}$ derived from the reconstruction formula \eqref{eqn:low-rank IB reconstruction}, we have 
\begin{align*}
{\rm II}^2_2\leq& \alpha^{-2}\|\big(\frac{|p|^2}{\epsilon^2}-1\big)\mathcal{F}^k(p;\gamma)\|^2_{L^2(B(0,\epsilon))}+\alpha^{-2}\|\big(\frac{|p|^2}{\epsilon^2}-1\big)\mathcal{F}^k(p;\eta)\|^2_{L^2(B(0,\epsilon))}\\
\leq& \alpha^{-2}\|\frac{|p|^2}{\epsilon^2}-1\|^2_{L^2(B(0,\epsilon))}\big(\sup_{|p|<\epsilon}|\mathcal{F}^k(p;\gamma)|^2+\sup_{|p|<\epsilon}|\mathcal{F}^k(p;\eta)|^2\big),
\end{align*}
where 
\begin{align*}
\|\frac{|p|^2}{\epsilon^2}-1\|^2_{L^2(B(0,\epsilon))} =\int_{B(0,\epsilon)}\Big|\frac{|p|^2}{\epsilon^2}-1\Big|^2 {\rm d} p =2\pi \int_{0}^{\epsilon}\Big|\frac{\rho^2}{\epsilon^2}-1\Big|^2 \rho{\rm d}\rho=\frac{\pi}{3}\epsilon^2
\end{align*}
and 
\begin{align*}
\sup_{|p|<\epsilon}|\mathcal{F}^k(p;\gamma)|^2+\sup_{|p|<\epsilon}|\mathcal{F}^k(p;\eta)|^2\leq& \big(\int_{B} |e^{i2kp\cdot y}| |\gamma(y)| {\rm d} y\big)^2+\big(\int_{B} |e^{i2kp\cdot y}| |\eta(y)| {\rm d} y\big)^2\\
\leq &|B|(\|\gamma\|^2_{L^2(B)}+\|\eta\|^2_{L^2(B)})=|B|\|\psi\|^2_{(L^2(B))^2}.
\end{align*}
Combing the above estimates and the identity $|B|=\pi$ gives
\begin{align*}
\|(I-\mathcal{K}_1K_1)\psi\|_{(L^2(B))^2}\leq 3^{-\frac12}\pi\alpha^{-1} \epsilon \|\psi\|_{(L^2(B))^2}+\delta_\alpha.
\end{align*}
\end{proof}

Recall that $\tilde{\psi}$ denotes the sum of the IBS and $$\mathcal{M}:=\max(\|\psi\|_{(L^2(B))^2},\|\tilde\psi\|_{(L^2(B))^2})\leq \mu_2^{-1} \big(1-\sqrt{1-(1+C_{\mathcal{K}_1})^{-1}}\big),$$ where 
$C_{\mathcal{K}_1}:=\nu_2\|\mathcal{K}_1\|_{(L^2(B))^2\to(L^2(B))^2}$, 
$\nu_{2}=2^\frac32|B| M^{-\frac12}$ and $\mu_{2} = 2^\frac32 M^{-\frac12}\big(\mu_0(k)+\mu_0(\ell k)\big)$, with $$\mu_0(k)=\max(1,k b(k), k^{-1}c(k),k^2 |B|^\frac12 a(k)),$$ $a(k)=\sup_{x\in B}\|G^k(x,\cdot)\|_{L^2(B)}\leq k^{-\frac12}$, $b(k)=\|F_1^k\|_{(L^2(B))^2\to L^2(B)}\leq |B|^\frac12 \big(3\sqrt{2}+\sqrt{\frac{8\pi}{3}}\big)k^\frac12$ and $c(k)=\|{\bf F_2^k}\|_{L^2(B)\to (L^2(B))^2}\leq |B|^\frac12 \big(3\sqrt{2}+\sqrt{\frac{8\pi}{3}}\big) k^\frac52$.   Based on the estimates in Lemmas \ref{lem:IBS_approx_err_Hoskins} and \ref{lem:K1_approx_err}, we can derive the following theorem on the approximation error of the IBS \eqref{eqn:IBS}.
\begin{theorem}\label{thm:IBS_approx_err}
Let the assumptions in Theorem \ref{thm:IBS_conv_2} hold and that the forward and IBS \eqref{eqn:IBS} converge. 
Then for any $\psi=(\gamma,\eta)^t\in (L^2(B))^2$ such that $\big(\|\gamma-\gamma^\alpha\|^2_{L^2(B)}+\|\eta-\eta^\alpha\|^2_{L^2(B)}\big)^\frac12\leq\delta_\alpha$ for some $\gamma^\alpha,\eta^\alpha\in \mbox{span}\{ \psi_{m,n,l}(\cdot;2k): |\alpha_{m,n}(2k)| \geq\alpha\}$ with $\psi_{m,n,l}$ being the disk PSWF discussed in Section \ref{sec: disk PSWFs}, 
the approximation error of the IBS \eqref{eqn:IBS} can be bounded by
\begin{align*}
\|\psi-\sum_{j=1}^{N} \mathcal{K}_j(\phi)\|_{(L^2(B))^2}\leq &2\mu_2 \big(\sqrt{16C_2^2+1}(1-C_{{\rm ratio}})\big)^{-1}C_{{\rm ratio}}^{N+1}\\
&+\Big(1+\big(1-(1-\mu_2\mathcal{M})^{-2}\big)C_{\mathcal{K}_1}\Big)^{-1}(3^{-\frac12}\pi\alpha^{-1} \epsilon \mathcal{M}+\delta_\alpha),
\end{align*}
where $C_2=\max(2,C_{\mathcal{K}_1})$ and $C_{{\rm ratio}}=\|\mathcal{K}_1\phi\|_{(L^2(B))^2}/r$.
\end{theorem}

\begin{proof}
It follows directly from Lemmas \ref{lem:IBS_approx_err_Hoskins} and \ref{lem:K1_approx_err}.
\end{proof}
\section{Numerical experiments}\label{sec:num}

In this section, we shall present preliminary numerical results to demonstrate the feasibility of the IBS to recover the perturbations $\eta$ and $\gamma$ from noisy scattering data (with $2\%$ relative noise) at two wave frequencies generated from the scattering problem \eqref{eqn:prob1}. 

\subsection{Data generation}
In the numerical experiments below, we consider a unit disk $\Omega=B$.
The perturbations of interest, i.e., $\gamma=a-1$ and $\eta=n-1$, are generated by 
\begin{enumerate}
\item [(i)]
\textbf{[Unseparated perturbations]} The weighted sums of Gaussian functions combined with a mask, given by 
\begin{equation}\label{eqn:dataset_s}
\gamma=\chi_\Omega \sum_{j=1}^{J} c_{a,j}e^{-\frac{(x - x_{a,j})^2 + (y - y_{a,j})^2}{2\sigma_{a,j}^2}} \quad \mbox{and} \quad \eta=\chi_\Omega\sum_{j=1}^{J} c_{n,j}e^{-\frac{(x - x_{n,j})^2 + (y - y_{n,j})^2}{2\sigma_{n,j}^2}},
\end{equation}
where $\chi_\Omega$ is the indicator function of $\Omega$, the peaks $$(x_{a,j},y_{a,j}) = (r_{a,j}\cos\theta_{a,j},r_{a,j}\sin\theta_{a,j})\quad \mbox{and}\quad(x_{n,j},y_{n,j}) = (r_{n,j}\cos\theta_{n,j},r_{n,j}\sin\theta_{n,j})$$ locate inside the circle of radius $0.5$ centered at the origin, and the standard deviations $\sigma_{a,j} = R_{a,j}(8\ln2)^{-\frac12}$ and $\sigma_{n,j} = R_{n,j}(8\ln2)^{-\frac12}$,
with $c_{a,j},c_{n,j}\in [0,J^{-1})$, $r_{a,j},r_{n,j}\in[0,0.5)$, $\theta_{a,j},\theta_{n,j}\in[0,2\pi)$, $R_{a,j}\in  \big(1-\max(|x_{a,j}|,|y_{a,j}|)\big)[0.3,1)$ and $R_{n,j}\in \big(1-\max(|x_{n,j}|,|y_{n,j}|)\big)[0.3,1)$ randomly sampled from uniform distributions. The bounds for $r_{a,j},r_{n,j}$ and $R_{a,j},R_{n,j}$ guarantee that the significant part of the Gaussian functions are within the unit circle and are not too concentrated in order to avoid small inhomogeneities, and the bounds for $c_{a,j},c_{n,j}$ ensure that the values of $a$ and $n$ are not too far away from $1$.
\item[(ii)]\textbf{[Separated perturbations]} 
The weighted sums of 2 Gaussian functions combined with a mask, given by 
\begin{equation}\label{eqn:dataset_pertu2}
\begin{array}{cc}
\gamma= &  \chi_\Omega  \big(c^{+}_{a}e^{-\frac{(x - x^+_{a})^2 + (y -y^+_{a})^2}{2(\sigma^+_{a})^2}}+c^{-}_{a}e^{-\frac{(x - x^-_{a})^2 + (y -y^-_{a})^2}{2(\sigma^-_{a})^2}}\big)\\
\eta= &  \chi_\Omega  \big(c^{+}_{n}e^{-\frac{(x - x^+_{n})^2 + (y -y^+_{n})^2}{2(\sigma^+_{n})^2}}+c^{-}_{n}e^{-\frac{(x - x^-_{n})^2 + (y -y^-_{n})^2}{2(\sigma^-_{n})^2}}\big),
\end{array}
\end{equation}
where $\chi_\Omega$ is the indicator function of $\Omega$, the peaks \begin{align*}
(x^{\pm}_{a},y^{\pm}_{a}) = &(r_{a}\cos\theta_{a} \pm 0.3,r_{a}\sin\theta_{a} \pm 0.3)\\
\mbox{and}\quad (x^{\pm}_{n},y^{\pm}_{n}) = &(r_{n}\cos\theta_{n}\pm 0.3,r_{n}\sin\theta_{n}\pm 0.3)
\end{align*}
locate inside the circle of radius $0.8$ centered at the origin, and the standard deviations $\sigma^{\pm}_{a} = R^{\pm}_{a}(8\ln2)^{-\frac12}$ and $\sigma^{\pm}_{n} = R^{\pm}_{n}(8\ln2)^{-\frac12}$,
with $c^{\pm}_{a},c^{\pm}_{n}\in [0,0.5)$, $r_{a},r_{n}\in[0,0.5)$, $\theta_{a},\theta_{n}\in[0,2\pi)$, $R^{\pm}_{a}\in  \big(1-\max(|x^{\pm}_{a}+0.2|,|y^{\pm}_{a}+0.2|)\big)[0.3,1)$ and $R^{\pm}_{n}\in \big(1-\max(|x^{\pm}_{n}+0.2|,|y^{\pm}_{n}+0.2|)\big)[0.3,1)$ randomly sampled from uniform distributions. 
The design of the peaks and the bounds for $R^{\pm}_{a}$ and $R^{\pm}_{n}$ help separate the significant part of the Gaussian functions, although not all of them are fully separated.
\end{enumerate}

The scattering data at wave frequency $k\in\{5,10,15\}$ are generated by solving the scattering problem \eqref{eqn:prob1}, where the radiation boundary conditions are implemented using the perfectly matched layer with order 4, using the Finite Element Method (implemented in NGsolve, available from https://github.com/NGSolve) with polynomials of order $4$ on a mesh of sizes $\frac{\pi}{4\sqrt{2}k}$ and $\frac{\pi}{4k}$ in $\Omega$ and in the air, respectively.

The exact perturbations $\gamma$ and $\eta$ are discretized in the computational domain $[-1,1)^2$ with a equispaced mesh of $N_{\rm out}\times N_{\rm out}$ points (i.e., $N_{\rm out}\times N_{\rm out}$ pixels), while
a single scattering data $u^\infty$ is first measured with $N_{\rm in}$ receivers and $N_{\rm in}$ sources with the same equiangular directions $\{\hat{x}_i,\hat{\theta}_j\}_{i,j=1}^{N_{\rm in}}$. 
Then the datum $u^\infty(\hat{x}_i,\hat{\theta}_j;k)$ on $\mathbb{S}\times\mathbb{S}$ is further transformed to an (approximately) equivalent datum evaluated at $p_{n}$, i.e., $u^\infty(p_n;k)$ on $B$, under the following transformation: 
Let $\{t_j,\omega_{t_j}\}_{j=0}^{T-1}$ be the set of Gauss-Legendre quadrature nodes and weights, and $\{\theta_i=\frac{2i\pi}{M},\omega_{\theta_i}=\frac{2\pi}{M}\}_{i=0}^{M-1}$ be the set of trapezoidal quadrature nodes and weights, then we identity $\big\{\sqrt{\frac{t_j+1}{2}} (\cos\theta_i,\sin \theta_i)^t \big\}_{j=0,~i=0}^{T-1,~M-1}$ by $\{p_n\}_{n=1}^{TM}$ 
and approximate $u(p_n)$ by 
$$
u(p_n) \approx u^\infty (\hat{x}_{i^*}, \hat{\theta}_{j^*};k), \quad \mbox{where}\quad
(i^*,j^*) =\mbox{argmin}_{i,j} \big\|p_n-\frac{\hat{\theta}_{j}-\hat{x}_{i}}{2} \big \|_2.$$
Note that the quadrature nodes satisfy that $|p_n|\neq 0$.
We refer the reader to \cite{ZhouAudibertMengZhang:2024} for the selection of $M$ and $T$. 
In the numerical experiments, we set $N_{\rm in}=N_{\rm out}=64$.

\subsection{Regularized inverse Born approximation}\label{sec:num_K1}

In the following, we provide several results on the capability of the regularized inverse Born approximation to recover sufficiently small perturbations of the form \eqref{eqn:dataset_s}, with $J=5$ and magnitudes less than $0.04$, from the scattering data $u^\infty$ at wave frequency $k\in\{5,10,15\}$ generated by solving the scattering problem \eqref{eqn:prob1}.
After generating, reorganizing and preprocessing the scattering data with $A(p_n)^{-1}$, we implement the reconstruction formula \eqref{eqn:low-rank IB reconstruction} with $c=2k$ and $\alpha=0.9\|\mathcal{F}^k\|_{L^2(B)\to L^2(B)}=0.9 |\alpha_{0,0}|$, which first projects $\{u(p_n)\}_{n=1}^{TM}$ onto the disk PSWFs; see \cite{ZhouAudibertMengZhang:2024} for the explicit algorithms and \cite{ZLWZ20,ZhouAudibertMengZhang:2024} for the evaluation of disk PSWFs.

\begin{figure}[htbp] 
\centering
Example 1:\qquad\qquad\qquad\qquad\qquad
Example 2: \\
\includegraphics[width=.39\linewidth]{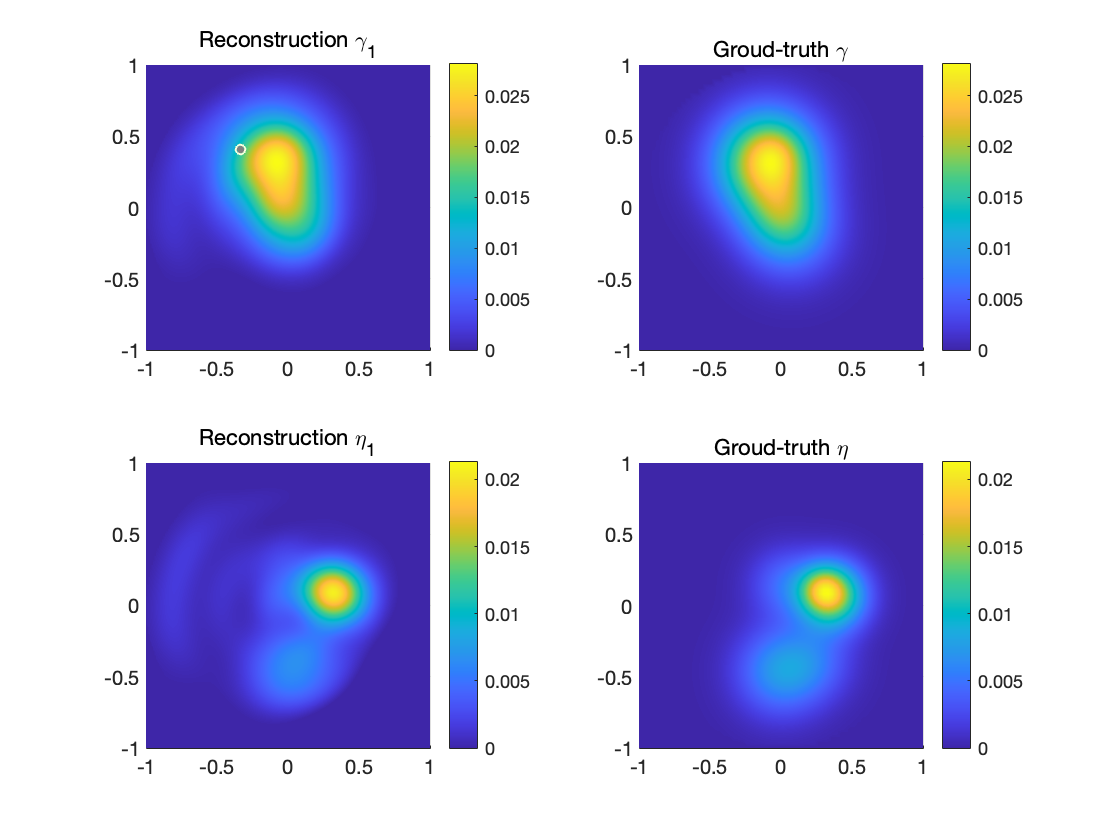} 
\includegraphics[width=.39\linewidth]{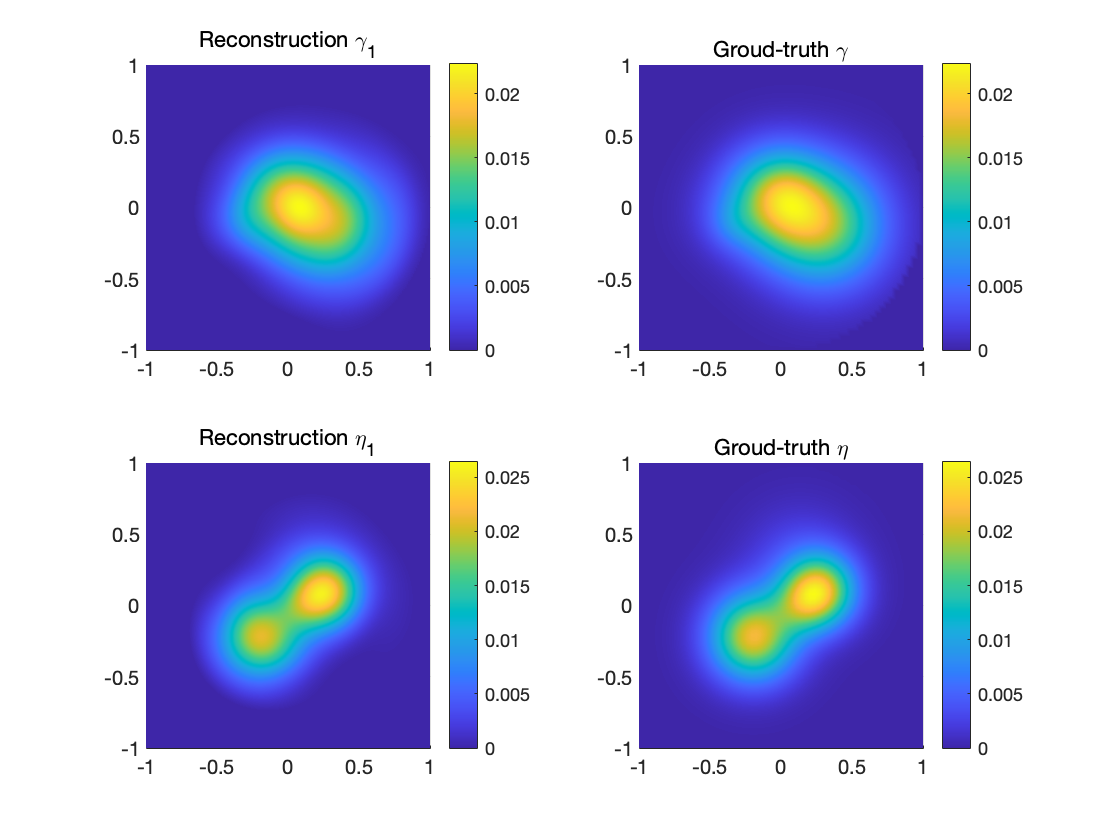} \\
Example 3:\qquad\qquad\qquad\qquad\qquad
Example 4: \\
\includegraphics[width=.39\linewidth]{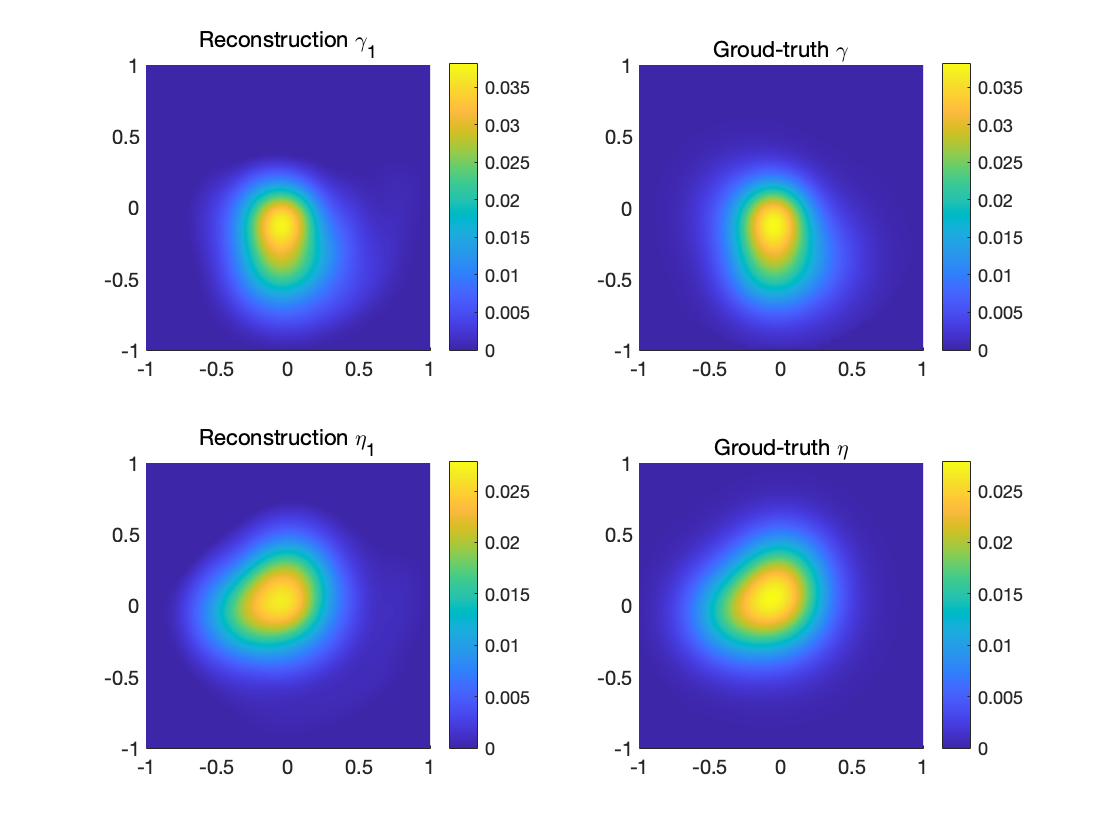} 
\includegraphics[width=.39\linewidth]{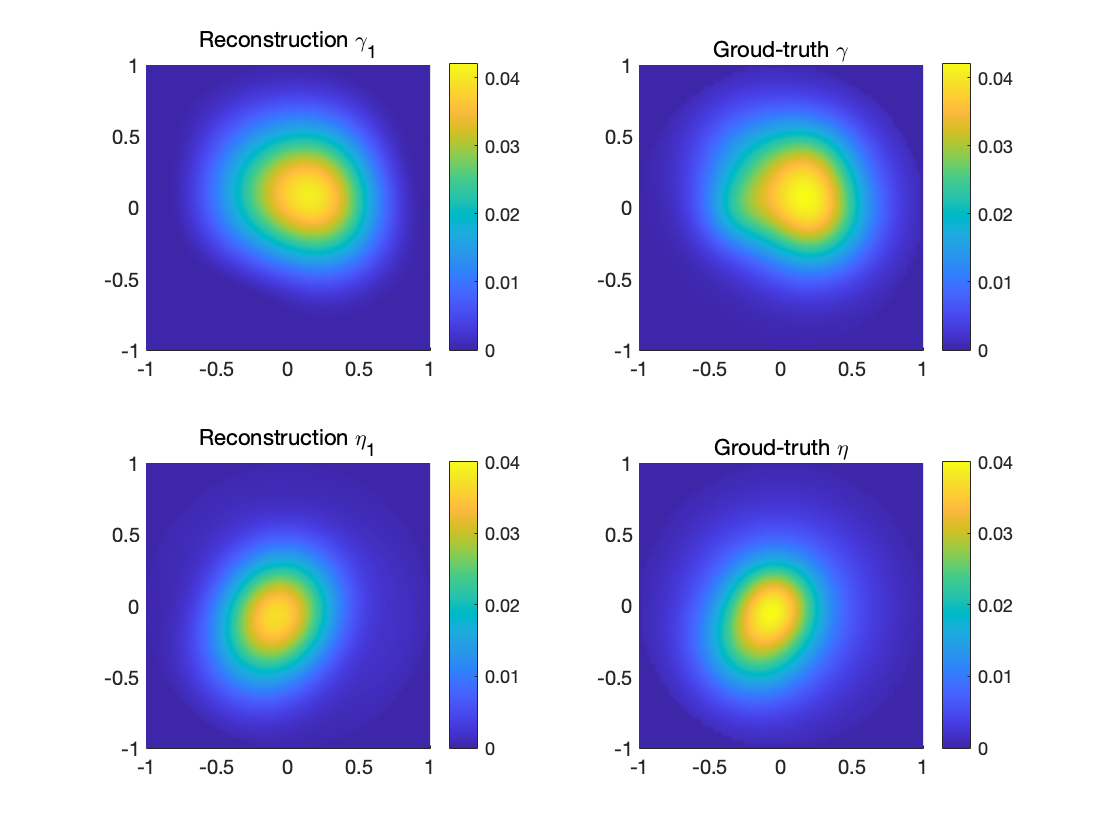} \\
\caption{Reconstructions of $\gamma$ and $\eta$ derived from the proposed regularized inverse Born approximation for different scattering data (with $2\%$ relative noise). }
\label{fig:born_low_rank}
\end{figure}

The reconstructions of the perturbations $\gamma$ and $\eta$ generated by the regularized inverse Born approximation are given in Figure \ref{fig:born_low_rank}. 
The results indicate that the proposed method has the capability to recover small perturbations with desirable quality. 
However, when the perturbations are relatively large, the inverse Born approximation which relies on linearized inverse problems  may fail to reconstruct the unknown inhomogeneities; see Figures \ref{fig:IBS} and \ref{fig:IBS_var_k_1}.

\subsection{Regularized inverse Born Series}\label{sec:num_IBS}
\begin{figure}[htbp] 
\centering
Example 1:\\
\includegraphics[width=0.99\linewidth]{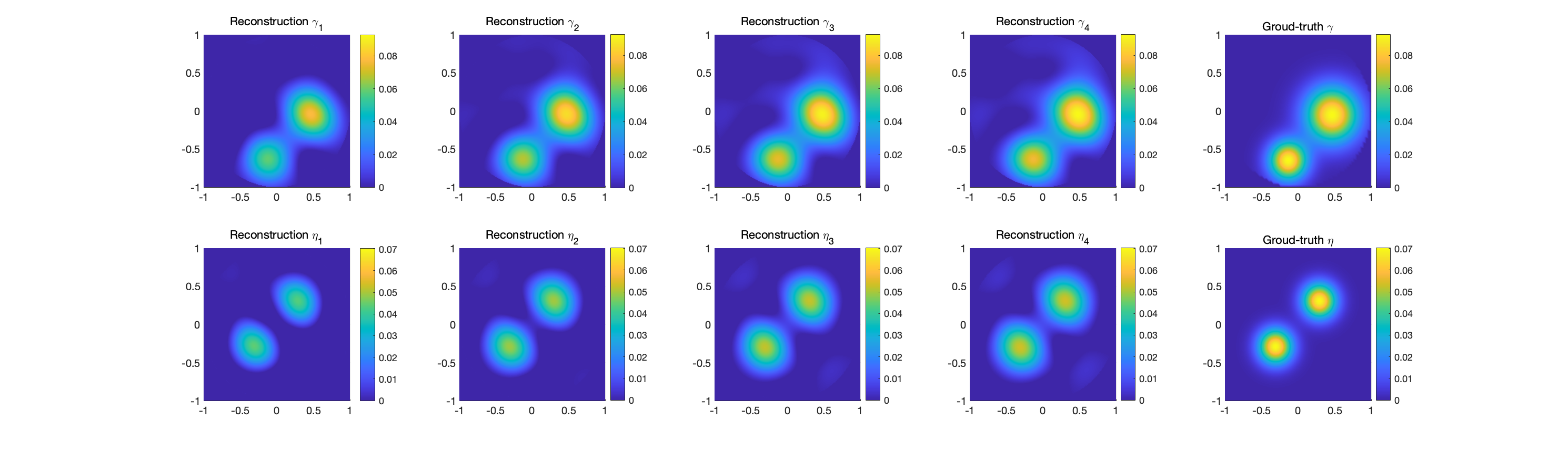} \\
Example 2:\\
\includegraphics[width=0.99\linewidth]{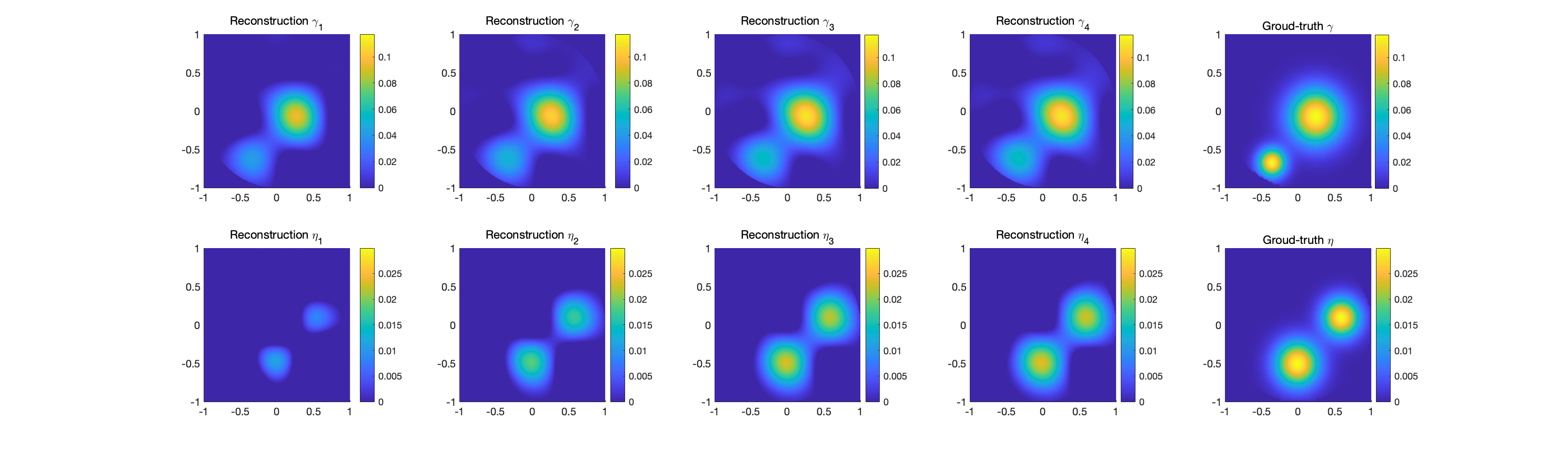} \\
Example 3:\\
\includegraphics[width=0.99\linewidth]{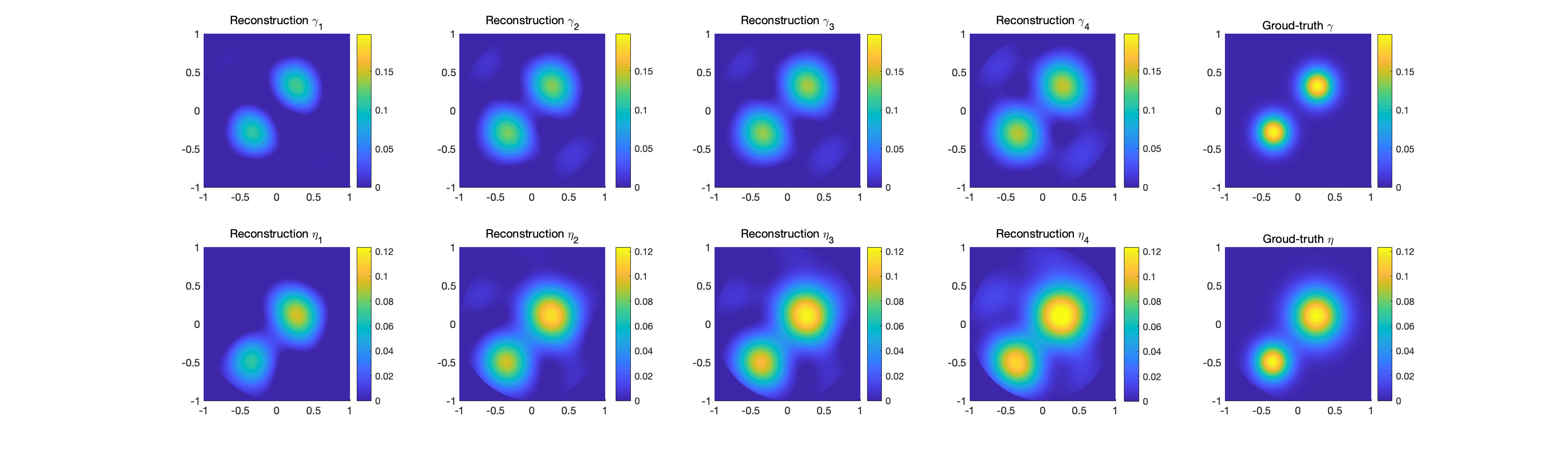} 
\caption{Reconstructions of $\gamma$ and $\eta$, i.e., $\gamma_i$ and $\eta_i$ that derived from the truncated IBS 
for different scattering data (with $2\%$ relative noise) at wave frequencies $k=5$ and $\ell k=10$, generated with the perturbations of the form \eqref{eqn:dataset_pertu2}.}
\label{fig:IBS_pertu2}
\end{figure}

\begin{figure}[htbp] 
\centering
Example 1:\\
\includegraphics[width=0.99\linewidth]{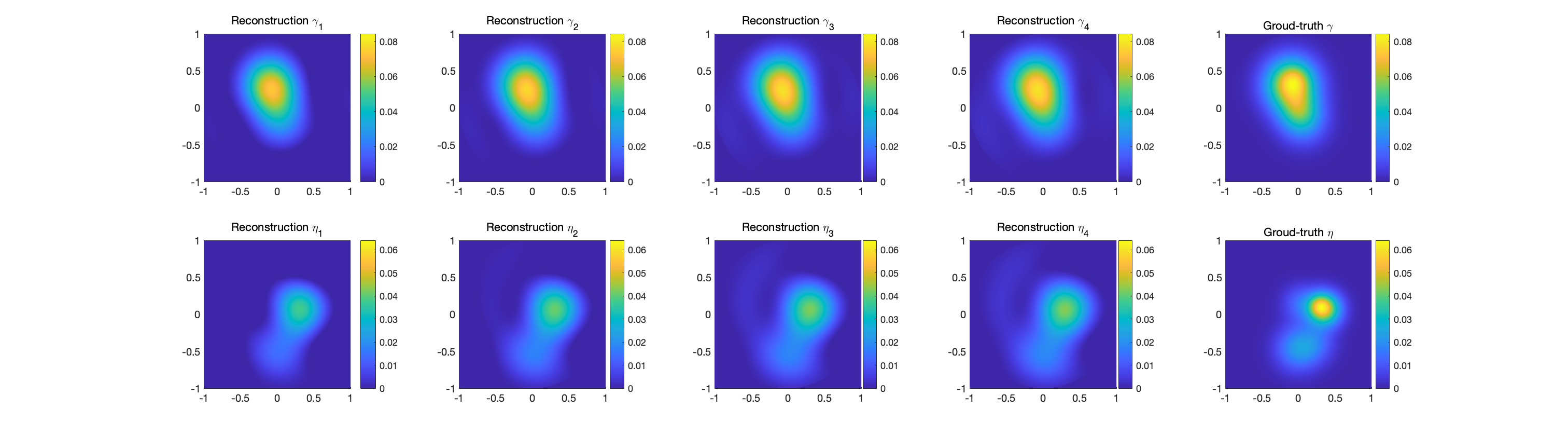} \\
Example 2:\\
\includegraphics[width=0.99\linewidth]{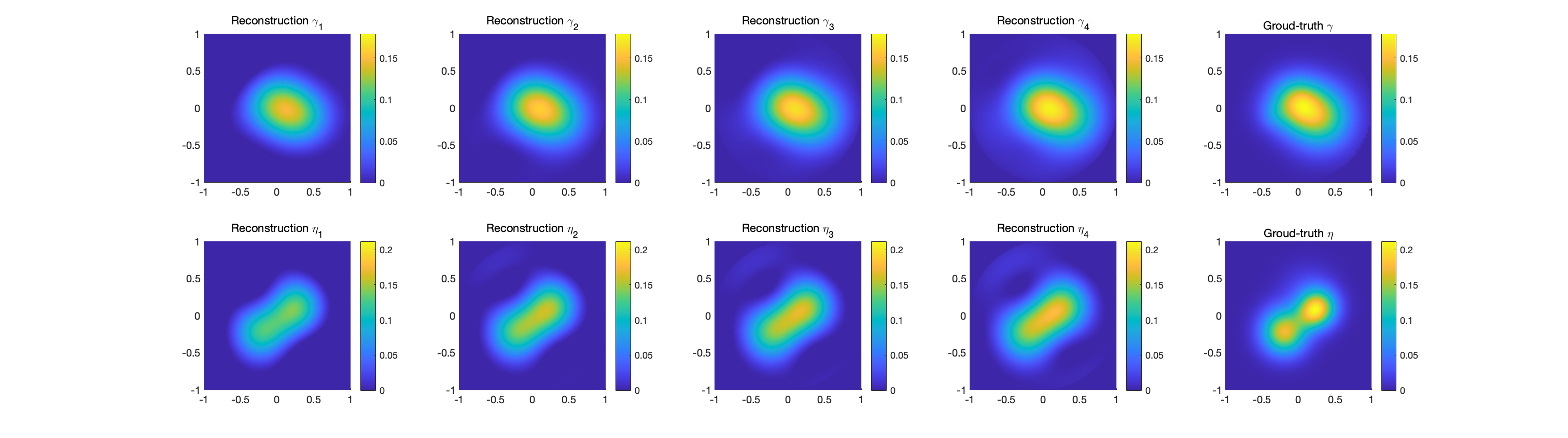} \\
Example 3:\\
\includegraphics[width=0.99\linewidth]{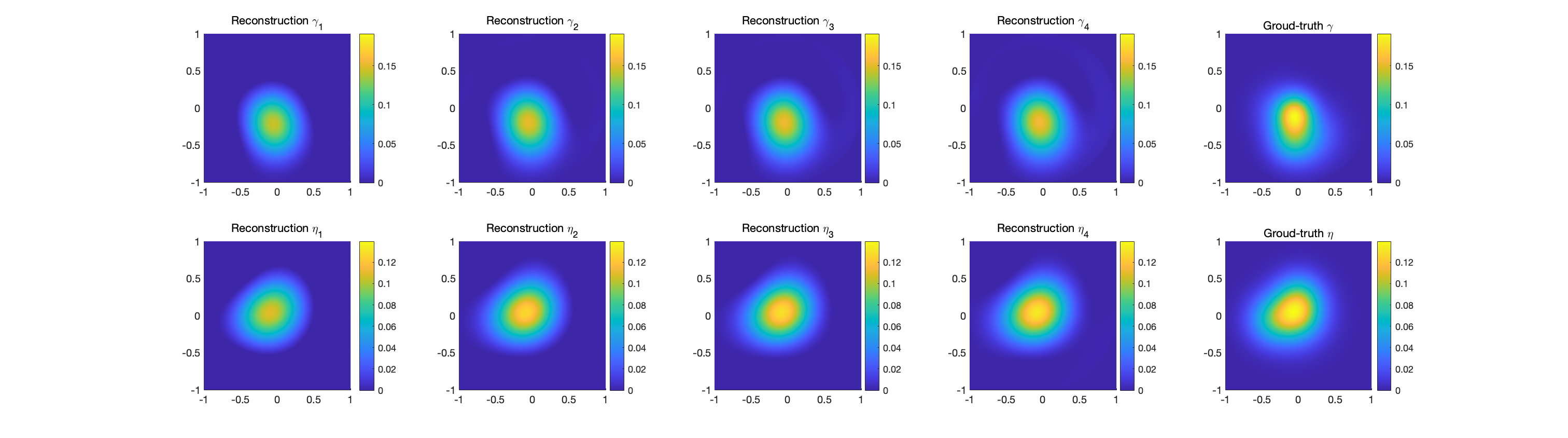} \\
Example 4:\\
\includegraphics[width=0.99\linewidth]{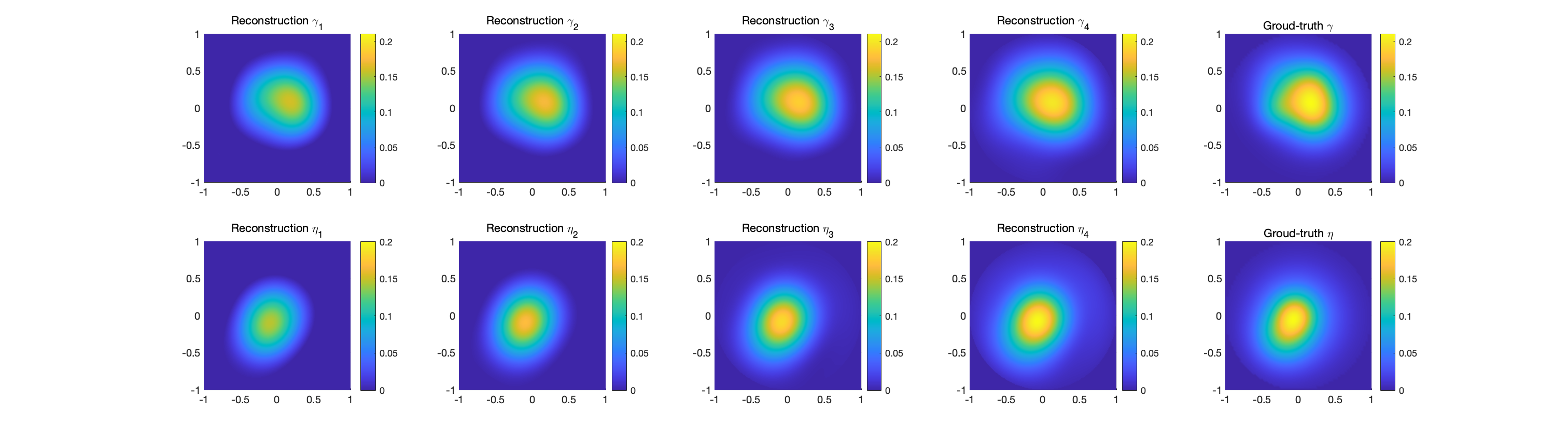} \\
\caption{Reconstructions of $\gamma$ and $\eta$, i.e., $\gamma_i$ and $\eta_i$ that derived from the truncated IBS 
for different scattering data (with $2\%$ relative noise) at wave frequencies $k=5$ and $\ell k=10$, generated with the perturbations of the form \eqref{eqn:dataset_s} ($J=5$).}
\label{fig:IBS}
\end{figure}

Now, we discuss the effect of IBS on reconstructing relatively large (unseparated and separated) perturbations of the forms \eqref{eqn:dataset_s} and \eqref{eqn:dataset_pertu2}, with $J=5$ and magnitudes between $0.08$ and $0.2$, from the scattering data $u^\infty$ at wave frequency $k\in\{5,10,15\}$ generated by solving the scattering problem \eqref{eqn:prob1}. 
Based on the estimates on the radius of convergence of IBS in Theorems \ref{thm:IBS_conv} and \ref{thm:IBS_conv_2}, we select the scattering data at $k=5$ and $\ell k=10$ to recover the perturbations $\gamma$ and $\eta$ for the convergence of IBS. 
The results of the reconstructions of separated $\gamma$ and $\eta$ are given in Figure \ref{fig:IBS}, while the results for unseparated $\gamma$ and $\eta$ are given in Figure \ref{fig:IBS_pertu2}. 
The reconstructions of $\gamma$ and $\eta$ derived from the truncated IBS $\sum_{j=1}^{i} \mathcal{K}_j(\phi)$ for different scattering data (with $2\%$ relative noise) are denoted by $\gamma_i$ and $\eta_i$, respectively. 

\begin{figure}[htbp] 
\centering
$\mathcal{K}_2$, $\mathcal{K}_3$ and $\mathcal{K}_4$ with $k=5$ and $\ell k=15$:\\
\includegraphics[width=0.99\linewidth]{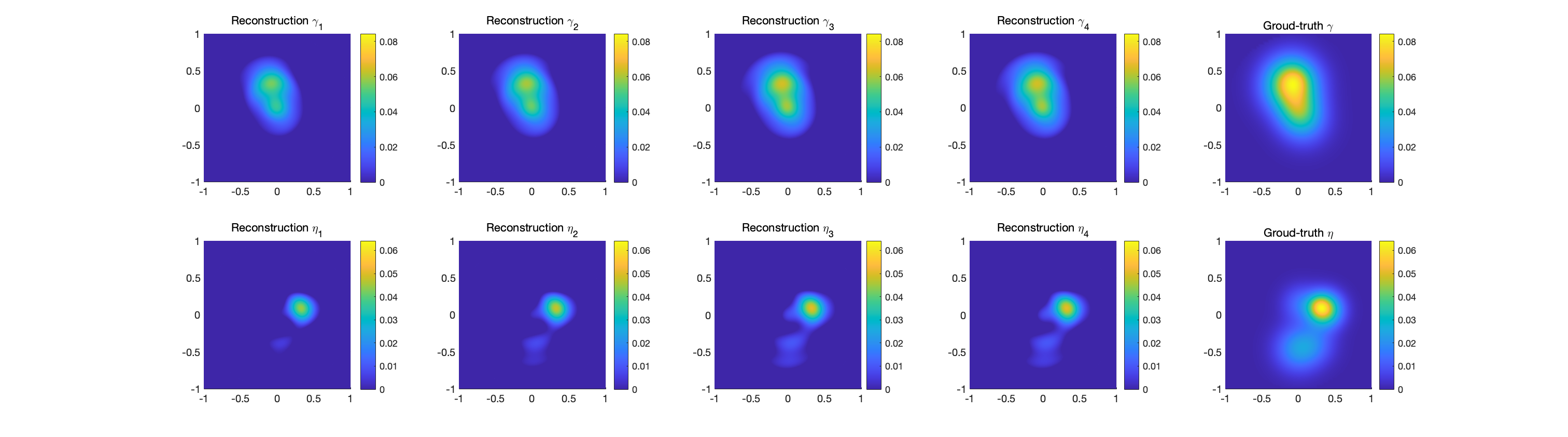} \\
$\mathcal{K}_2$ and $\mathcal{K}_3$ with $k=10$ and $\ell k=15$ \& $\mathcal{K}_4$ with $k=5$ and $\ell k=15$:\\
\includegraphics[width=0.99\linewidth]{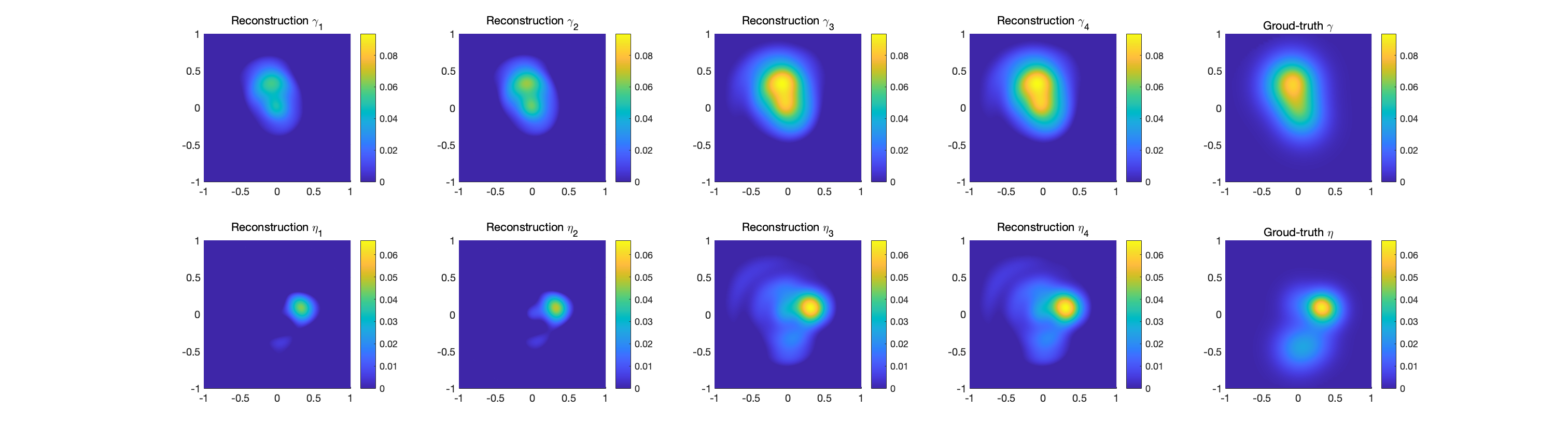} \\
\caption{Reconstructions of $\gamma$ and $\eta$, i.e., $\gamma_i$ and $\eta_i$ that derived from the truncated IBS 
for different scattering data (with $2\%$ relative noise) at wave frequencies $k=10$ and $\ell k=15$, generated with the perturbations of Example 1 in Figure \ref{fig:IBS}.
}
\label{fig:IBS_var_k_1}
\end{figure}

\begin{figure}[htbp] 
\centering
$\mathcal{K}_2$, $\mathcal{K}_3$ and $\mathcal{K}_4$ with $k=5$ and $\ell k=15$:\\
\includegraphics[width=0.99\linewidth]{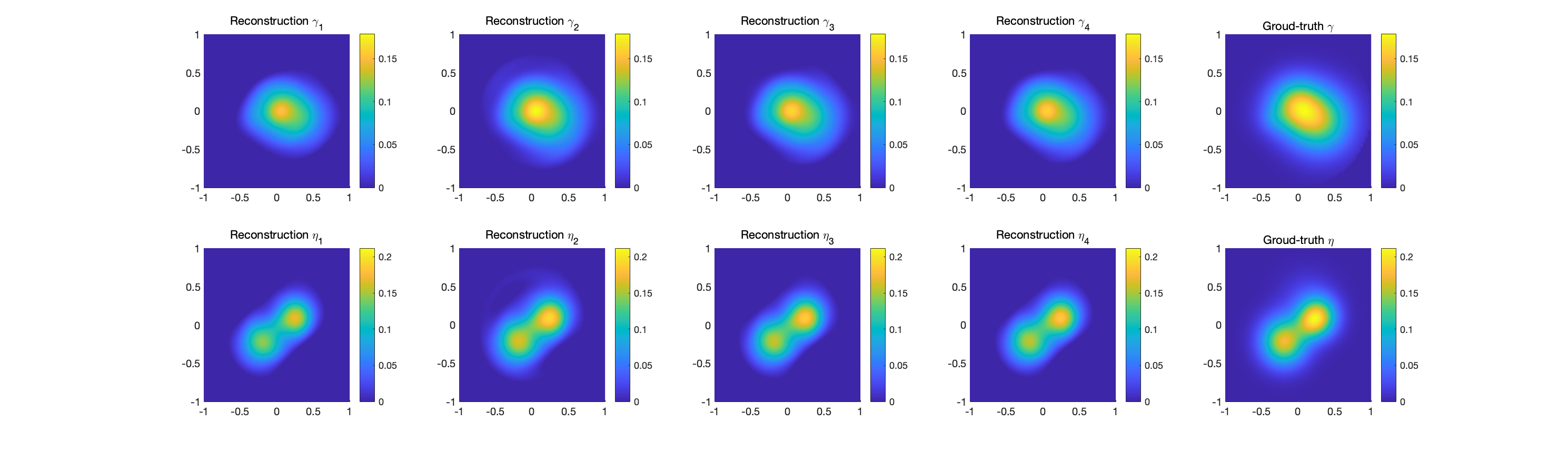} \\
$\mathcal{K}_2$ with $k=5$ and $\ell k=15$ \& $\mathcal{K}_3$ and $\mathcal{K}_4$ with $k=5$ and $\ell k=10$:\\
\includegraphics[width=0.99\linewidth]{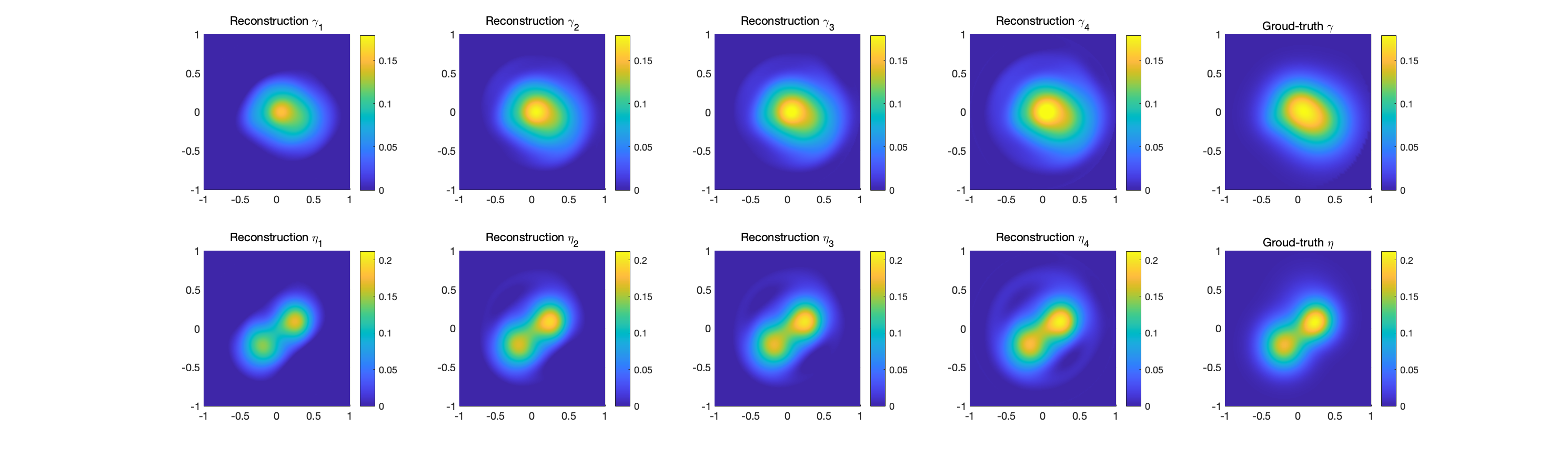} \\
\caption{Reconstructions of $\gamma$ and $\eta$, i.e., $\gamma_i$ and $\eta_i$ that derived from the truncated IBS 
for different scattering data (with $2\%$ relative noise) at wave frequencies $k=10$ and $\ell k=15$, generated with the perturbations of Example 2 in Figure \ref{fig:IBS}.
}
\label{fig:IBS_var_k_2}
\end{figure}

\begin{figure}[htbp] 
\centering
$\mathcal{K}_2$, $\mathcal{K}_3$ and $\mathcal{K}_4$ with $k=5$ and $\ell k=15$:\\
\includegraphics[width=0.99\linewidth]{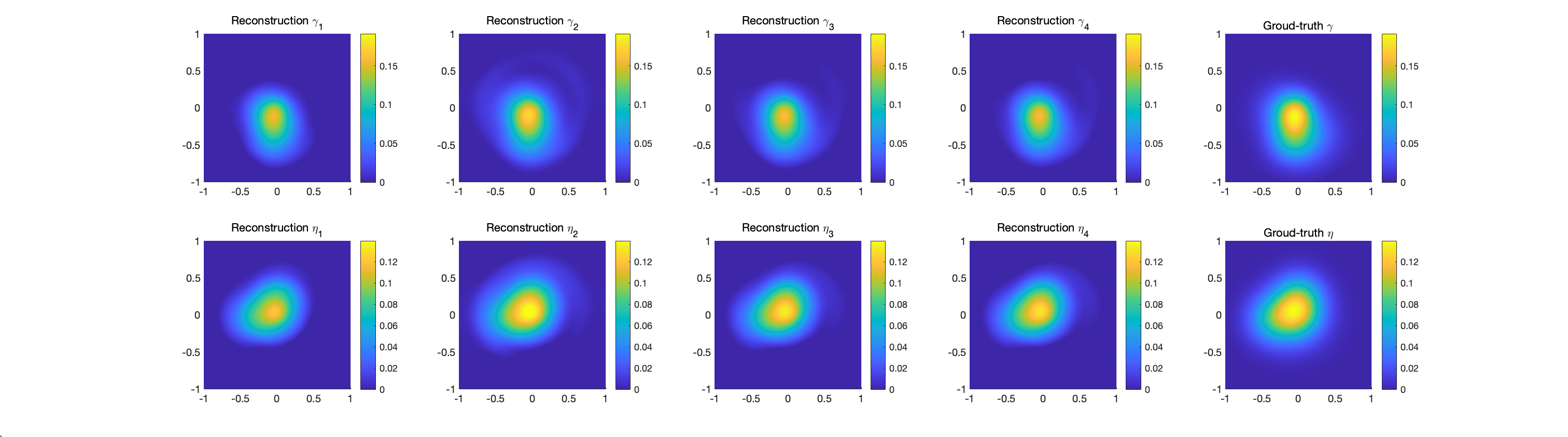} \\
$\mathcal{K}_2$ with $k=5$ and $\ell k=15$ \& $\mathcal{K}_3$ and $\mathcal{K}_4$ with $k=5$ and $\ell k=10$:\\
\includegraphics[width=0.99\linewidth]{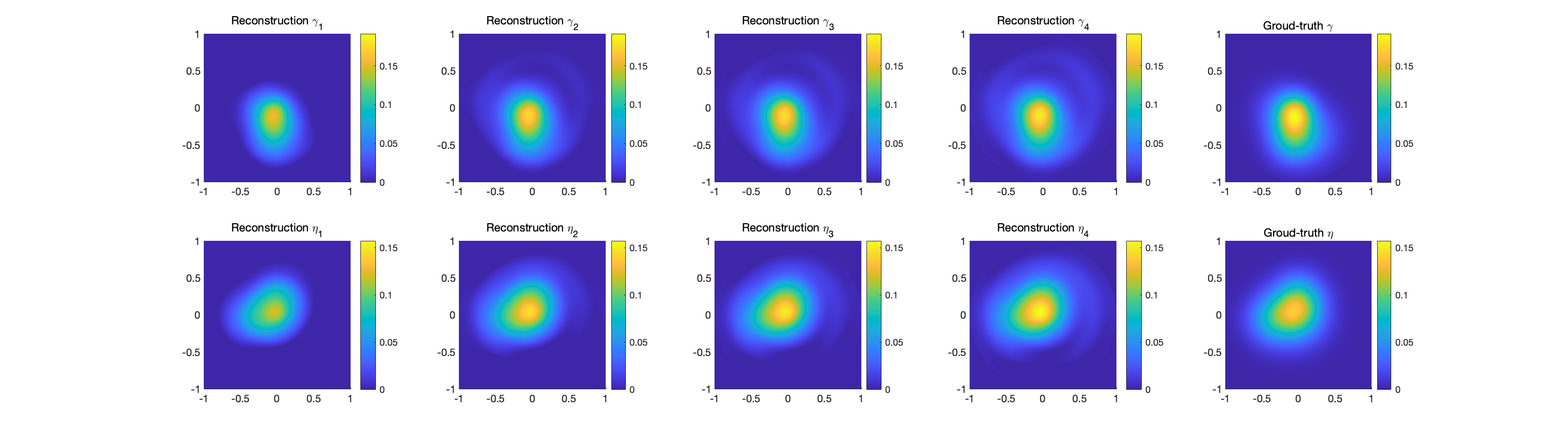} \\
\caption{Reconstructions of $\gamma$ and $\eta$, i.e., $\gamma_i$ and $\eta_i$ that derived from the truncated IBS 
for different scattering data (with $2\%$ relative noise) at wave frequencies $k=10$ and $\ell k=15$, generated with the perturbations of Example 3 in Figure \ref{fig:IBS}.
}
\label{fig:IBS_var_k_3}
\end{figure}

An improvement (either less or more) in the reconstructions is observed in these examples. Specifically, the inverse Born approximation locates the perturbations, while the subsequent terms in the series help recover the magnitude of the perturbations.
The quality of the reconstructions from IBS closely depends on that of the inverse Born approximation.
Based on the observation discussed in \cite{ZhouAudibertMengZhang:2024}, the reconstruction obtained from the proposed regularized inverse Born approximation becomes better when the frequency $k$ increases. 
As observed in Examples 1--3 of Figure \ref{fig:IBS}, the current wave frequencies may not be high enough for accurate reconstructions.
To improve the quality of the reconstructions from IBS, we use the scattering data at higher wave frequencies in the regularized inverse Born approximation. 
However, higher wave frequencies may lead to the divergence of the IBS due to the estimates in Theorems \ref{thm:IBS_conv} and \ref{thm:IBS_conv_2}.
Thus, we investigate the reconstructions of IBS from the scattering data at wave frequencies $k=10$ and $\ell k=15$, where the first-order term, i.e., the regularized inverse Born approximation, is evaluated at $k=10$ and $\ell k=15$ while other terms are evaluated at $k=5$ or $k=10$ and $\ell k=10$ or $\ell k=15$. 
The results for Examples 1, 2 and 3 in Figure \ref{fig:IBS} are given in Figures \ref{fig:IBS_var_k_1}, \ref{fig:IBS_var_k_2} and \ref{fig:IBS_var_k_3}, respectively.
Note that although we compute $\{\mathcal{K}_j(\phi)\}_{j=2}^i$ at different wave frequencies, the input data $\phi$ depends only on the scattering data at $k=10$ and $\ell k=15$.

On the basis of the above observation, the wave frequencies for the IBS are important for both its convergence and the quality of the reconstructions. 
A proper balance is necessary for accurate recovery of the inhomogeneities.


\appendix
\section{Auxiliary estimates}\label{app:estimate}
In this appendix, we give several auxiliary estimates that have been used in the convergence analysis of the IBS \eqref{eqn:IBS}.
First, we state the proof of Proposition \ref{prop:nu_mu}.
\begin{proof}[\textbf{Proof of Proposition \ref{prop:nu_mu}}]
For any $j\geq 1$, it follows directly from a slightly modified version of recursion \eqref{eqn:us_j} for $u^s_j$, where $\eta$ and $\gamma$ are replaced by $\eta_j$ and $\gamma_j$ respectively, as well as from the definitions of the integral operators $F_0^k$ and $F_1^k$ in \eqref{eqn:F0} and \eqref{eqn:F1}, that
\begin{align*}
&\|u^s_{j}\|_{L^2(B)}=\|F_1^k(\gamma_{j}\nabla u^s_{j-1})+F_0^k(\eta_{j} u^s_{j-1})\|_{L^2(B)}\\
\leq & \|F_1^k\|_{(L^2(B))^2\to L^2(B)}\|\gamma_{j}\nabla u^s_{j-1}\|_{(L^2(B))^2}+|B|^\frac12\|F_0^k(\eta_{j} u^s_{j-1})\|_{L^\infty(B)}\\
\leq &\|F_1^k\|_{(L^2(B))^2\to L^2(B)}\|\gamma_{j}\|_{L^{\infty}(B)}\|\nabla u^s_{j-1}\|_{(L^2(B))^2}+|B|^\frac12 k^2\sup_{x\in B}|\int_{B} G^k(x,y)\eta_{j}(y)u^s_{j-1}(y) {\rm d}y|\\
\leq & b(k)\|\gamma_{j}\|_{L^{\infty}(B)}\|\nabla u^s_{j-1}\|_{(L^2(B))^2} +k^2 |B|^\frac12 a(k)\|\eta_{j}\|_{L^{\infty}(B)}\|u^s_{j-1}\|_{L^2(B)}
\end{align*}
where $a(k)=\sup_{x\in B}\|G^k(x,\cdot)\|_{L^2(B)}$ and $b(k)=\|F_1^k\|_{(L^2(B))^2\to L^2(B)}<\infty$. 
Indeed, $F_1^k: (L^2(B))^2\to L^2(B)$ has a weakly singular kernel and hence is a compact operator with operator norm $\|F_1^k\|_{(L^2(B))^2\to L^2(B)}<\infty$.
Similarly, with the modified version of recursion \eqref{eqn:us_j} for $\nabla u^s_j$, the definitions of ${\bf F_2^k}$ and ${\bf F_3^k}$ in \eqref{eqn:F2} and \eqref{eqn:F3}, and the estimate $\|{\bf F_3^k}\|_{(L^2(B))^2\to (L^2(B))^2}\leq 1$ stated in \cite[Lemma 2.2]{ArridgeMoskowSchotland:2012}, there holds
\begin{align*}
\|\nabla u^s_{j}\|_{(L^2(B))^2}=&\|{\bf F_3^k}(\gamma_{j} \nabla u^s_{j-1})+{\bf F_2^k}(\eta_{j} u^s_{j-1})\|_{(L^2(B))^2}\\
\leq & \|\gamma_{j} \nabla u^s_{j-1}\|_{(L^2(B))^2}+\|{\bf F_2^k}\|_{L^2(B)\to (L^2(B))^2}\|\eta_{j}\|_{L^\infty(B)}\|u^s_{j-1}\|_{L^2(B)}\\
\leq & \|\gamma_{j}\|_{L^\infty(B)} \|\nabla u^s_{j-1}\|_{(L^2(B))^2}+c(k)\|\eta_{j}\|_{L^\infty(B)}\|u^s_{j-1}\|_{L^2(B)},
\end{align*}
where $c(k)=\|{\bf F_2^k}\|_{L^2(B)\to (L^2(B))^2}<\infty$.
We define $a_j(k):=\|u^s_{j}\|_{L^2(B)}$ and $b_j(k):=\|\nabla u^s_{j}\|_{(L^2(B))^2}$, and obtain the following coupled recursions:
\begin{align}\label{eqn:coupled us_j}
\begin{aligned}
a_j(k)\leq& b(k)\|\gamma_{j}\|_{L^{\infty}(B)}b_{j-1}(k)+k^2|B|^\frac12 a(k)\|\eta_{j}\|_{L^{\infty}(B)}a_{j-1}(k),\\
b_j(k)\leq& \|\gamma_{j}\|_{L^\infty(B)}b_{j-1}(k)+ c(k)\|\eta_{j}\|_{L^\infty(B)}a_{j-1}(k)   
\end{aligned}
\end{align}
with the initial terms $a_0(k)=\|u^i\|_{L^2(B)}=|B|^\frac12$ and $b_0(k)=\|\nabla u^i\|_{L^2(B)}\leq k |B|^\frac12$.

Next, by the recursion \eqref{eqn:u_j} for $u_j(p;k)=u_j(\hat{x},\hat{\theta};k)$ with the relations $q-p=\hat{x}\in\mathbb{S}$,
we derive that
\begin{align*}
\|u_j(p;k)\|_{L^2(B)}\leq &|B|^\frac12 \|ik^{-1}\mathcal{F}^{p,k}_0(\gamma_{j} (q-p)\cdot\nabla u^s_{j-1})+\mathcal{F}^{p,k}_0(\eta_{j} u^s_{j-1})\|_{L^\infty(B)}\\
\leq &|B|^\frac12 \big( k^{-1}\|\gamma_{j}\|_{L^\infty(B)}b_{j-1}(k)+ \|\eta_{j}\|_{L^{\infty}(B)} a_{j-1}(k)\big).
\end{align*}
We suppress the subscripts of $\|\eta_{j}\|_{L^{\infty}(B)}$ and $\|\gamma_{j}\|_{L^{\infty}(B)}$ in the following analysis. By noting that $\tilde{\phi}_j=(u_j(p;k),u_j(\ell^{-1}p;\ell k)\big)^t$, the above inequality and the coupled recursions \eqref{eqn:coupled us_j} imply that
\begin{align*}
&\|\tilde{\phi}_{j}\|_{(L^2(B))^2}=(\|u_{j}(p;k)\|_{L^2(B)}^2+\|u_{j}(\ell^{-1}p;\ell k)\|_{L^2(B)}^2)^\frac12\\
\leq & \|u_{j}(p;k)\|_{L^2(B)}+\|u_{j}(\ell^{-1}p;\ell k)\|_{L^2(B)}\\
\leq &|B|^\frac12\big((k^{-1} \|\gamma_{j}\|  b_{j-1}(k)+  \|\eta_{j}\|  a_{j-1}(k))+((\ell k)^{-1} \|\gamma_{j}\|  b_{j-1}(\ell k)+  \|\eta_{j}\|  a_{j-1}(\ell k))\big)\\
:=&|B|^\frac12\big({\rm I}_{j}(k)+{\rm I}_{j}(\ell k)\big),
\end{align*}
where
\begin{align*}
{\rm I}_{j}(k)=&k^{-1} \|\gamma_{j}\|  b_{j-1}(k)+  \|\eta_{j}\|  a_{j-1}(k)\\
\leq& \big(\|\gamma_{j}\| +  kb(k)\|\eta_{j}\|\big) k^{-1}\|\gamma_{j-1}\|  b_{j-2}(k)\\
&+ \big(k^{-1} c(k)\|\gamma_{j}\|+k^2|B|^\frac12 a(k) \|\eta_{j}\| \big)\|\eta_{j-1}\|  a_{j-2}(k)\\
\leq & \max(1,k b(k),k^{-1}c(k),k^2 |B|^\frac12 a(k))(\|\gamma_{j}\|+\|\eta_{j}\|) {\rm I}_{j-1}(k)\\
:= & \mu_0(k)(\|\gamma_{j}\|+\|\eta_{j}\|) {\rm I}_{j-1}(k),\\
{\rm I}_j(\ell k)\leq &  \mu_0(\ell k)(\|\gamma_{j}\|+\|\eta_{j}\|) {\rm I}_{j-1}(\ell k)
\end{align*}
with $\mu_0(k)=\max(1,k b(k), k^{-1}c(k),k^2 |B|^\frac12 a(k))$.
Let $\tilde{\psi}_j=(\gamma_j,\eta_j)^t$ for any $j\geq 1$, the Cauchy–Schwarz inequality gives $\|\gamma_{j}\|+\|\eta_{j}\|\leq \sqrt{2}(\|\gamma_{j}\|^2+\|\eta_{j}\|^2)^\frac12=\sqrt{2}\|\tilde{\psi}_j\|_{(L^\infty(B))^2}$ and hence
\begin{align*}
{\rm I}_{j}(k)\leq& \sqrt{2}\mu_0(k)\|\tilde{\psi}_{j}\|_{(L^\infty(B))^2}{\rm I}_{j-1}(k) \leq (\sqrt{2}\mu_0(k))^{j-1} \Pi_{i=2}^{j} \|\tilde{\psi}_i\|_{(L^\infty(B))^2} {\rm I}_1(k)\\
\leq & (\sqrt{2}\mu_0(k))^{j-1} \Pi_{i=2}^{j} \|\tilde{\psi}_i\|_{(L^\infty(B))^2} (k^{-1} \|\gamma_1\|  b_{0}(k)+  \|\eta_1\|  a_{0}(k))\\
\leq & |B|^\frac12(\sqrt{2}\mu_0(k))^{j-1} \Pi_{i=2}^{j} \|\tilde{\psi}_i\|_{(L^\infty(B))^2} (\|\gamma_1\|+  \|\eta_1\|)\\
\leq & \sqrt{2}|B|^\frac12(\sqrt{2}\mu_0(k))^{j-1} \Pi_{i=1}^{j} \|\tilde{\psi}_i\|_{(L^\infty(B))^2}.
\end{align*}
Finally, by combining above estimates, we obtain
\begin{align*}
\|\tilde{\phi}_{j}\|_{(L^2(B))^2}\leq & |B|^\frac12 ({\rm I}_{j-1}(k)+{\rm I}_{j-1}(\ell k))\\
\leq & |B|\sqrt{2}^{j} \big( \mu_0(k)^{j-1} + \mu_0(\ell k)^{j-1} \big) \Pi_{i=1}^{j} \|\tilde{\psi}_i\|_{(L^\infty(B))^2}\\
\leq & |B| \sqrt{2}^{j} \big(\mu_0(k)+\mu_0(\ell k)\big)^{j-1} \Pi_{i=1}^{j} \|\tilde{\psi}_i\|_{(L^\infty(B))^2}\\
:=& \nu_{\infty} \mu_{\infty}^{j-1}\Pi_{i=1}^{j} \|\tilde{\psi}_i\|_{(L^\infty(B))^2},
\end{align*}
where $\nu_{\infty}=\sqrt{2}|B|$ and $\mu_{\infty} = \sqrt{2} \big(\mu_0(k)+\mu_0(\ell k)\big)$.
This completes the proof.
\end{proof}

Next, we provide upper bounds of the coefficients $a(k)$, $b(k)$ and $c(k)$ within Proposition \ref{prop:nu_mu} in the following lemma.
\begin{lemma}\label{lem:abc}
Let $a(k)=\sup_{x\in B}\|G^k(x,\cdot)\|_{L^2(B)}$, $b(k)=\|F_1^k\|_{(L^2(B))^2\to L^2(B)}$ and $c(k)=\|{\bf F_2^k}\|_{L^2(B)\to (L^2(B))^2}$, then for any $k>\frac12$, there hold
\begin{align*}
a(k)\leq& \frac{\sqrt{2k+1}}{2k}\leq k^{-\frac12},\quad b(k)\leq |B|^\frac12 \big(\frac{3}{2k^{\frac32}} +\sqrt{\frac{8\pi}{3}}\big)k^{\frac12}
\leq |B|^\frac12 \big(3\sqrt{2}+\sqrt{\frac{8\pi}{3}}\big)k^\frac12\\
\mbox{and} \quad c(k)\leq&  |B|^\frac12 \big(\frac{3}{2k^{\frac32}} +\sqrt{\frac{8\pi}{3}}\big)k^\frac52\leq |B|^\frac12 \big(3\sqrt{2}+\sqrt{\frac{8\pi}{3}}\big) k^\frac52.
\end{align*}
\end{lemma}
\begin{proof}
First, we evaluate $a(k)$. 
By noting the upper bound of $|H_0^{(1)}(\rho)|$ that 
\begin{align*}
|H_0^{(1)}(\rho)|\leq 
\left\{
\begin{aligned}1-\ln \rho,&\quad 0<\rho\leq 1,\\
(\frac{2}{\pi \rho})^\frac12,&\quad 1< \rho,\\
\end{aligned}\right.
\end{align*} 
we derive from the definition $G^k(x,y)=\frac{i}{4} H_0^{(1)}(k |x-y|)$ that 
\begin{align*}
a(k)^2=&\sup_{x\in B}\|G^k(x,\cdot)\|^2_{L^2(B)}\leq \int_{2B} |G^k(0,y)|^2{\rm d}y
\leq \frac1{16} \int_{2B} |H_0^{(1)}(k|y|)|^2{\rm d}y\\
\leq & \frac1{16} \int_{0}^{2\pi}\int_{0}^{2} |H_0^{(1)}(k\rho)|^2 \rho{\rm d}\rho {\rm d}\theta
\leq \frac{\pi}{8} \int_{0}^{2} |H_0^{(1)}(k\rho)|^2 \rho{\rm d}\rho
\\
\leq&  \frac{\pi}{8} (\int_{0}^{k^{-1}}+\int_{k^{-1}}^2) |H_0^{(1)}(k\rho)|^2 \rho{\rm d}\rho 
\leq \frac{\pi}{8} \big(\int_{0}^{k^{-1}} (1-\ln (k \rho))^2\rho{\rm d}\rho +\int_{k^{-1}}^2 \frac{2}{\pi k}{\rm d}\rho\big)\\ 
\leq& \frac{\pi}{8} \big(k^{-2}\int_{0}^{1} (1-\ln \rho)^2\rho{\rm d}\rho +\frac{2(2-k^{-1})}{\pi k}\big)
\leq \frac{\pi}{8} \big(\frac{5}{4 k^2}+\frac{2(2k-1)}{\pi k^2}\big)
\leq \frac{2k+1}{4k^2},
\end{align*}
which implies $a(k)\leq \frac{\sqrt{2k+1}}{2k}$.

\vspace{2mm}
Next, we consider $b(k)=\|F_1^k\|_{(L^2(B))^2\to L^2(B)}$. 
By the definition of $F_1^k$ in \eqref{eqn:F1} and the identity $\nabla_x G^k(x,y)=-\frac{ik(x-y)}{4|x-y|} H_1^{(1)}(k|x-y|)$ , there holds
\begin{align*}
\|F_1^k\|_{(L^2(B))^2\to L^2(B)}=&\sup_{\|{\bf f}\|_{(L^2(B))^2}=1}\|F_1^k({\bf f})\|_{L^2(B)}
\leq |B|^\frac12 \sup_{\|{\bf f}\|_{(L^2(B))^2}=1} \|F_1^k({\bf f})\|_{L^\infty(B)}\\
\leq &|B|^\frac12 \sup_{\|{\bf f}\|_{(L^2(B))^2}=1} \sup_{x\in B} \int_{B}|\nabla_{x} G^k(x,y)| |{\bf f}(y)|\; {\rm d}y\\
\leq &\frac{k}{4}|B|^\frac12 \sup_{\|{\bf f}\|_{(L^2(B))^2}=1}\sup_{x\in B}\int_{B}|H_1^{(1)}(k|x-y|)||{\bf f}(y)|\; {\rm d}y\\
:=&\frac{k}{4}|B|^\frac12 \sup_{\|{\bf f}\|_{(L^2(B))^2}=1}\sup_{x\in B}{\rm I}(x).
\end{align*}
Then the upper bound of $|H_1^{(1)}(\rho)|$ that 
\begin{align*}
|H_1^{(1)}(\rho)|\leq 
\left\{
\begin{aligned}\frac{3}{\pi \rho}\quad,&\quad 0<\rho\leq 1,\\
(\frac{3}{\pi \rho})^\frac12,&\quad 1< \rho,\\
\end{aligned}\right.
\end{align*} 
indicates
\begin{align*}
{\rm I}(x)= &\int_{B}|H_1^{(1)}(k|x-y|)||{\bf f}(y)|\; {\rm d}y\\
\leq &\frac{3}{\pi k}\int_{B}\chi_{\{y:|x-y|\leq k^{-1}\}}|x-y|^{-1} |{\bf f}(y)|\; {\rm d}y+\sqrt{\frac{3}{\pi k}}\int_{B}|x-y|^{-\frac12}|{\bf f}(y)|\; {\rm d}y\\
:= &\frac{3}{\pi k}\int_{B}h_1(x,y) |{\bf f}(y)|\; {\rm d}y+\sqrt{\frac{3}{\pi k}}\int_{B}h_2(x,y)|{\bf f}(y)|\; {\rm d}y\\
:=&\frac{3}{\pi k} H_1(x;|{\bf f}(y)|)+\sqrt{\frac{3}{\pi k}}H_2(x;|{\bf f}(y)|).
\end{align*}
Next, for any $x\in B$, the Schur's test ensure that 
\begin{align*}
\|H_1\|_{L^2(B)\to L^2(B)}\leq& \sup_{x\in B}\int_{B} h_1(x,y){\rm d}y
\leq \int_{2B} \chi_{\{y:|y|\leq k^{-1}\}}|y|^{-1}{\rm d}y\\
\leq& \int_{0}^{2\pi}\int_0^{k^{-1}} \rho^{-1}\rho {\rm d}\rho {\rm d} \theta\leq 2\pi k^{-1}
\end{align*}
and 
\begin{align*}
\|H_2\|_{L^2(B)\to L^2(B)}\leq& \sup_{x\in B}\int_{B} h_2(x,y){\rm d}y
\leq \int_{2B}|y|^{-\frac12}{\rm d}y\\
\leq& \int_{0}^{2\pi}\int_{0}^2 \rho^{-\frac12}\rho {\rm d}\rho {\rm d} \theta\leq \frac{8\sqrt{2}\pi}{3}.
\end{align*}
Finally, by combining the above estimates and the identity $\|{\bf f}\|_{(L^2(B))^2}=\||{\bf f}|\|_{L^2(B)}$, we obtain
\begin{align*}
b(k)=&\|F_1^k\|_{(L^2(B))^2\to L^2(B)}\leq \frac{k}{4}|B|^\frac12 \sup_{\||{\bf f}|\|_{L^2(B)}=1}\sup_{x\in B}\big(\frac{3}{\pi k} H_1(x;|{\bf f}(y)|)+\sqrt{\frac{3}{\pi k}}H_2(x;|{\bf f}(y)|)\big)\\
\leq &\frac{k}{4}|B|^\frac12 \big(\frac{3}{\pi k} \|H_1\|_{L^2(B)\to L^2(B)}+\sqrt{\frac{3}{\pi k}}\|H_2\|_{L^2(B)\to L^2(B)}\big)\\
\leq &|B|^\frac12 \big(\frac{3}{2k^{\frac32}} +\sqrt{\frac{8\pi}{3}}\big)k^{\frac12}.
\end{align*}
Following the analysis for $b(k)$, we derive from the definition of ${\bf F_2^k}$ in \eqref{eqn:F2} that 
\begin{align*}
c(k)=\|{\bf F_2^k}\|_{L^2(B)\to (L^2(B))^2}\leq  |B|^\frac12 \big(\frac{3}{2k^{\frac32}} +\sqrt{\frac{8\pi}{3}}\big)k^{\frac52}.
\end{align*}
This completes the proof.
\end{proof}

\bibliographystyle{abbrv}
\bibliography{inv_born}

\end{document}